\documentclass[11pt]{article}
\usepackage[utf8]{inputenc}

\usepackage{amsthm, amsfonts, amssymb, amsmath,enumitem, natbib, bbm, mathabx}
\usepackage{cases, amsthm}
\usepackage{fullpage}
\usepackage{mathtools}
\usepackage{mathrsfs} 
\usepackage{wasysym} 
\usepackage{verbatim}
\usepackage{fancyhdr}
\usepackage{graphicx}
\usepackage{bbm}
\usepackage{caption}
\usepackage{subcaption}
\usepackage{algorithm}
\usepackage{algorithmic}
\usepackage{authblk}
\usepackage{setspace}
\let\Algorithm\algorithm
\renewcommand\algorithm[1][]{\Algorithm[#1]\setstretch{1.6}}
\usepackage[colorlinks=true, pdfstartview=FitV,
linkcolor=blue, citecolor=blue, urlcolor=blue]{hyperref}
\usepackage{cleveref}

\newtheorem{definition}{Definition}[section]
\newtheorem{theorem}[definition]{Theorem}
\newtheorem{proposition}[definition]{Proposition}
\newtheorem{lemma}[definition]{Lemma}
\newtheorem{corollary}[definition]{Corollary}
\newtheorem{example}{Example}[section]

\theoremstyle{remark}
\newtheorem*{remark}{Remark}
\newtheorem*{claim}{Claim}

\usepackage{accents}

\newcommand{\medbias}{\mathrm{Med}\mbox{-}\mathrm{bias}}

\DeclareMathOperator*{\argmax}{arg\,max}
\DeclareMathOperator*{\argmin}{arg\,min}

\newcommand{\sign}{\mathrm{sgn}}

\newcommand{\norm}[1]{\left\lVert#1\right\rVert}

\usepackage{xcolor}
\newcommand{\E}{\mathbbm{E}}

\newcommand{\g}{\gamma}

\def\ed { \stackrel{d}{=} }

\def\R{\mathbb{R}}

\def\E{\mathbb{E}}
 
\def\Z{\mathbb{Z}}
\def\d{\mathbb{D}}
\def\D{\mathbb{D}}
\def\g{\mathbb{G}}
\def\X{\mathbb{X}}

\DeclareMathOperator*{\median}{median \;}

\def\ed { \stackrel{d}{=} }

\def\R{\mathbb{R}}
\def\P{\mathbb{P}}

\def\E{\mathbb{E}}

\usepackage[top=1in, bottom=1in, left=1in, right=1in]{geometry}
\title{On the Symmetry of Limiting Distribution of M-estimators}

\author[1]{Arunav Bhowmick}
\author[2]{Arun Kumar Kuchibhotla}
\affil[1]{Indian Statistical Institute, Kolkata}
\affil[2]{Department of Statistics \& Data Science, Carnegie Mellon University}
\date{}

\begin{document}
\maketitle

\begin{abstract}  
Many functionals of interest in statistics and machine learning can be written as minimizers of expected loss functions. Such functionals are called $M$-estimands, and can be estimated by $M$-estimators --- minimizers of empirical average losses. Traditionally, statistical inference (e.g., hypothesis tests and confidence sets) for $M$-estimands is obtained by proving asymptotic normality of $M$-estimators centered at the target. However, asymptotic normality is only one of several possible limiting distributions and (asymptotically) valid inference becomes significantly difficult with non-normal limits. In this paper, we provide conditions for the symmetry of three general classes of limiting distributions, enabling inference using HulC~\citep{kuchibhotla2024hulc}.
\end{abstract}

\section{Introduction} \label{intro}

In many estimation problems, the functional of interest, $\theta_0 = \theta(P)$, is defined as the minimizer of the expected value of a suitable loss function, $L(\cdot, W)$, where $W$ is a random variable with the underlying distribution $P\in\mathcal{P}$, for a collection of data distributions $\mathcal{P}$. Formally,
\begin{equation} \label{define_theta}
    \theta_0 \equiv \theta(P) := \argmin_{\phi \in \Theta} \E_P [L(\phi, W)],
\end{equation}
where $\Theta$ is some parameter space; we allow the parameter space to be a subset of any metric space. We implicitly assume that there is a unique minimizer $\theta(P)$ for the optimization problem~\eqref{define_theta}. Traditional inference methods such as Wald interval, bootstrap, and subsampling rely on an estimator $\widetilde{\theta}_n$ based on $n$ IID observations $W_1, \ldots, W_n$ from $P$ which when centered at $\theta_0$ and properly normalized converges in distribution, as $n\to\infty$. Often, but not always, $\widetilde{\theta}_n$ is chosen to be the empirical minimizer (i.e., $M$-estimator). Assuming the parameter space $\Theta$ is known, an $M$-estimator is defined as 
\begin{equation} \label{define_theta_hat}
    \widehat{\theta}_n \in \argmin_{\phi \in \Theta} \frac{1}{n} \sum_{i=1}^{n}  L(\phi, W_i),
\end{equation}
where $\{W_i: 1\le i\le n\}$ is a sequence of IID random variables from $P$. If the minimizer is unique in~\eqref{define_theta_hat}, then $\widehat{\theta}_n$ is well-defined and if there are multiple global minimizers, then $\widehat{\theta}_n$ can be chosen at random. If the minimum is not attained in~\eqref{define_theta_hat}, then one can choose $\widehat{\theta}_n$ to be any point which is an $\varepsilon_n$-minimizer of $\phi\mapsto n^{-1}\sum_{i=1}^n L(\phi, W_i)$~\citep[Page 1303]{hess1996epi}. If, for all $P\in\mathcal{P}$, $r_n(\widehat{\theta}_n - \theta(P))$ converges in distribution to a random variable $W_P$, then (1) Wald intervals rely on a parametric form for the distribution of $W_P$ to derive valid inference, and (2) bootstrap, or more generally, resampling methods estimate the limiting distribution $W_P$ to derive valid inference. These methods tend to not yield {\em uniformly} valid inference (uniform over $P\in\mathcal{P}$) without strong continuity assumptions on the law of $W_P$ with respect to $P\in\mathcal{P}$~\citep[Sec. 8.5]{van2000asymptotic}. Such continuity assumptions often fail in non-standard or irregular problems (exemplified in the following).

Although normal limiting distribution commonly arises in parametric and semiparametric statistics, it is far from being the only possible limiting distribution. In case of mean zero normal limiting distributions, continuity of $W_P$ with respect to $P$ can be judged by continuity of the variance parameter. For the large classes of limiting distributions we will discuss, no such simple method exists.   Derivation of limiting distribution of $M$-estimator is not an easy task, in general. 

To provide valid inference in statistical problems where the continuity conditions could fail, \cite{kuchibhotla2024hulc} devised a new method of constructing confidence intervals called HulC. For univariate functionals of interest $\theta_0$ and $\alpha\in[0, 1]$, the $(1-\alpha)$ HulC confidence interval is $\widehat{\mathrm{CI}}_{n,\alpha}^{\mathrm{HulC}} := [\min_{1\le j\le B_{\alpha}}\widehat{\theta}^{(j)},\,\max_{1\le j\le B_{\alpha}}\widehat{\theta}^{(j)}]$ where $B_{\alpha} = \lceil\log_2(2/\alpha)\rceil$ and $\widehat{\theta}^{(j)}, 1\le j\le B_{\alpha}$ are independent estimators of $\theta_0.$ The main result of~\cite{kuchibhotla2024hulc} proves that HulC yields asymptotically valid confidence intervals if the median bias of the estimators $\widehat{\theta}^{(j)}$ converges to zero as $n\to\infty$, where the median bias is defined as
\begin{equation}\label{eq:median-bias}
\medbias_{P}(\widehat{\theta}_n; \theta(P)) = \left(\frac{1}{2} - \min_{s\in\{-1, 1\}}\mathbb{P}_P(s(\widehat{\theta}_n - \theta(P)) \ge 0)\right)_+.
\end{equation}
If $\medbias_P(\cdot;\cdot)$ converges uniformly over all $P\in\mathcal{P}$, HulC confidence interval is uniformly valid. This can be more easily satisfied than continuity of limiting distribution with respect to $P\in\mathcal{P}$; see~\cite{kuchibhotla2024hulc,kuchibhotla2023median} for examples and details. Note that if $r_n(\widehat{\theta}_n - \theta(P))\overset{d}{\to} W_P$ for $P\in\mathcal{P}$, then $\medbias_P(\widehat{\theta}_n; \theta(P))$ converges to $(1/2 - \min\{\mathbb{P}(sW_P \ge 0): s\in\{-1, 1\}\})_+,$ if $\mathbb{P}(W_P = 0) = 0$; see Appendix~\ref{appsec:proof-of-median-unbiasedness}. A multivariate generalization of HulC is discussed in~\cite{jain2023rectangular}. Median bias and its multivariate generalizations capture the idea that the estimator should not prefer one direction of $\theta_0$ from another direction; in the univariate case, over-estimation and under-estimation are equally likely. Hence, symmetry of the limiting distribution is a sufficient, but not a necessary condition for limiting median bias to be zero. 

The main goal of this paper is to identify sufficient conditions for symmetry of limiting distribution of $M$-estimators and investigate the extent to which such conditions are necessary. Even though the median bias of an estimator can converge to zero without converging in distribution~\citep[Page 1501]{pfanzagl1970asymptotic}, we restrict ourselves to cases where the $M$-estimator when properly normalized converges in distribution. We study the symmetry of three main classes of limiting distributions developed in~\cite{pflug1995asymptotic} and discuss statistical estimators fitting into these classes.

The remaining article is organized as follows. In Section~\ref{sec:limiting-distribution}, we discuss a general technique for obtaining the limiting distribution of $M$-estimators, provide examples, and also present required background for the asymptotics of constrained $M$-estimators; this discussion follows the works of~\cite{pflug1995asymptotic},~\cite{geyer1994asymptotics},~\cite{shapiro2000asymptotics}, and~\cite{knight1999epi}. In Section~\ref{sec:sufficiency}, we provide our first main result discussing a simple sufficient condition for the symmetry of min36imizers of stochastic processes and verify this condition in some examples. In Section~\ref{sec:necessity}, we study the necessity of the stated sufficient conditions for the symmetry of one class of limiting distributions. Finally, we conclude the article with some remarks in Section~\ref{sec:conclusions}. 

Our sufficiency results pave way for asymptotic valid inference via HulC for the wide range of (irregular/non-standard) problems. (As will be clear, this large class of problems contain the well-known cube-root estimators.) Our necessity results show the cases where HulC (as it is defined) is provably invalid.
\section{Limiting Distributions of $M$-estimators}\label{sec:limiting-distribution}
\subsection{Derivation of the Limiting Distribution}
Recall the definitions of $M$-estimand and $M$-estimator as
\[
\theta(P) = \argmin_{\phi\in\Theta}\,\mathbb{E}_P[L(\phi, W)],\quad\mbox{and}\quad \widehat{\theta}_n = \argmin_{\phi \in\Theta}\,\frac{1}{n}\sum_{i=1}^n L(\phi, W_i).
\]
For notational convenience, set $\mathbb{M}(\phi) = \mathbb{E}_P[L(\phi, W)]$ and $\mathbb{M}_n(\phi) = n^{-1}\sum_{i=1}^n L(\phi, W_i)$.  
For the following discussion, we assume $\Theta$ is a fixed non-stochastic set as $n$ changes and that $\Theta\subseteq\Omega$ for a location-scale invariant space $\Omega$; for a discussion on the limiting distributions for the cases with data-dependent or sample size dependent constraints, see~\cite{knight1999epi}. If the loss function $L(\cdot, \cdot)$ is differentiable in the first argument and minumum is attained in the interior of $\Theta$, then the estimand and the estimator can be written as zeros of some (normal) equations. This can then by combined with Taylor series expansion to derive asymptotic normality of the estimator. When the loss function is non-differentiable or when the minimum is attained at the boundary of $\Theta$ such a proof technique does not work. The general strategy can be understood as follows.

For any two non-negative sequences $\{r_n\}, \{\rho_n\}$, we have
\begin{equation} \label{define_Tn}
\begin{split}
    r_n(\widehat{\theta}_n - \theta_0) &= \argmin_{u \in r_n(\Theta - \theta_0)}\, \rho_n (\mathbb{M}_n(\theta_0 + u/r_n) - \mathbb{M}_n(\theta_0))\\
    &=  \argmin_{u\in\Omega}\, \mathbb{D}_n(u) + \mathbb{S}_n(u) + \mathbb{X}_n(u),
\end{split}
\end{equation}
where 
\begin{align*}
    \mathbb{D}_n(u) ~&:=~ \rho_n(\mathbb{M}(\theta_0 + u/r_n) - \mathbb{M}(\theta_0))\quad (\mbox{deterministic part}),\\
    \mathbb{S}_n(u) ~&:=~ \rho_n((\mathbb{M}_n - \mathbb{M})(\theta_0 + u/r_n) - (\mathbb{M}_n - \mathbb{M})(\theta_0)),\quad(\mbox{stochastic part})\\
    \mathbb{X}_n(u) ~&:= \begin{cases}0, &\mbox{if }u\in r_n(\Theta - \theta_0),\\
    +\infty, &\mbox{otherwise,}\end{cases} \quad(\mbox{constraint part)}.
\end{align*}
The function $\mathbb{X}_n(\cdot)$ is sometimes referred to as ``optimization theory indicator'' function~\citep[Example 2.67]{bonnans2013perturbation}. From this point onward, we will refer to such functions as indicator functions. Assuming convergence of $\mathbb{D}_n(\cdot), \mathbb{X}_n(\cdot)$ to functions $\mathbb{D}(\cdot), \mathbb{X}(\cdot)$ and convergence in distribution of $\mathbb{S}_n(\cdot)$ to $\mathbb{S}(\cdot)$ in an appropriate sense to ensure convergence of minimizers, we get 
\begin{equation}\label{eq:limiting-distribution}
r_n(\widehat{\theta}_n - \theta_0) ~\overset{d}{\to}~ \argmin_{u\in\Omega} \mathbb{D}(u) + \mathbb{S}(u) + \mathbb{X}(u).
\end{equation}
The choice of $r_n$ and $\rho_n$ should be made so as to ensure convergence of $\mathbb{D}_n(\cdot)$ and $\mathbb{S}_n(\cdot)$ to ``tight'' processes as $n\to\infty$. Also, the type of convergence of random processes to ensure convergence of minimizers should be chosen appropriately. Usually, one would use uniform convergence on compacts as in~\cite{van1996weak}, but for constrained problems, this is not suitable and one should use epi-convergence~\citep{rockafeller1998variational,knight1999epi}. Epi-convergence is considered the weakest form of convergence of functions that imply convergence of minimizers; see~\citet[page 2]{burke2017epi} and~\citet[Chapter 7.E]{rockafeller1998variational}. A sequence of non-random functions $\{f_n\}_{n\ge0}$ is said to epi-converge to a function $f$ if
\begin{align*}
    \liminf_{n\to\infty}\, f_n(x_n) ~&\ge~ f(x)\quad\mbox{for every sequence }x_n\to x,\\
    \limsup_{n\to\infty}\, f_n(x_n) ~&\le~ f(x)\quad\mbox{for some sequence }x_n\to x.
\end{align*}
Epi-convergence can be equivalently defined as convergence of epi-graphs of $f_n$ to that of $f$~\citep[Proposition 7.2]{rockafeller1998variational}. The relation between pointwise, uniform, and epi-convergence of a sequence of functions is rather delicate as detailed in~\citet[Chapter 7]{rockafeller1998variational} but under relatively mild conditions, epi-convergence implies pointwise convergence and is implied by uniform convergence. Also, see~\cite{kall1986approximation}. Moreover, Theorem 7.33 of~\cite{rockafeller1998variational} proves that epi-convergence under mild conditions implies $\argmin_x f_n(x) \subseteq \argmin_{x}\, f(x)$ without requiring the set of minimizers to be singletons. For random functions~\cite{salinetti1986convergence} defined epi-convergence in distribution for the purpose of concluding~\eqref{eq:limiting-distribution}; also, see~\cite{knight1999epi} for a review.

The limit function $\mathbb{X}(\cdot)$ of $\mathbb{X}_n$ is by design a function that only takes values $0$ and $+\infty$. Hence, $\mathbb{X}(\cdot)$ is the indicator function of some set $\mathcal{T}_{\Theta}(\theta_0)$, which is a tangent cone of $\Theta$ at $\theta_0$. We now briefly describe different types of tangent cones, following Chapter 4 of~\cite{aubin2009set}; also, see~\cite{geyer1994asymptotics} and~\cite{shapiro2000asymptotics}. For any $x\in\Omega$ and $S\subseteq\Omega$, $\mbox{dist}(x, S) = \inf_{z\in S}\|x - z\|$ for some norm $\|\cdot\|$. For any collection of sets $S(t)\subseteq\Omega$ indexed by $t\in\mathbb{R}^k$, define the limsup and liminf sets as
\begin{align*}
    \limsup_{t\to t_0}\,S(t) &:= \left\{x\in\Omega:\, \liminf_{t\to t_0}\,\mbox{dist}(x, S(t)) = 0\right\},\\
    \liminf_{t\to t_0}\,S(t) &:= \left\{x\in\Omega:\, \limsup_{t\to t_0}\,\mbox{dist}(x, S(t)) = 0\right\}.
\end{align*}
With this background in place, different tangent cones are defined as follows. (It is worth noting that the notation of tangent cones is not always consistent in the notation and we follow the notation of~\cite{aubin2009set}.) The {\em contingent} (Bouligand) tangent cone of $\Theta$ at $\theta_0$ is defined
\[
T_{\Theta}(\theta_0) ~:=~ \limsup_{t\uparrow\infty}\,t(\Theta - \theta_0).
\]
An alternative description is that $\vartheta\in T_{\Theta}(\theta_0)$ if and only if there exists a sequence $\{k_m\}$ diverging to $\infty$ and a sequence $\{\theta_m\}$ in $\Theta$ converging to $\theta_0$ such that $k_m(\theta_m - \theta_0) \to \vartheta$.
The {\em intermediate} or {\em adjacent} or {\em inner} tangent cone of $\Theta$ at $\theta_0$ is defined as
\[
{T}_{\Theta}^{\flat}(\theta_0) ~:=~ \liminf_{t\uparrow\infty}\, t(\Theta - \theta_0).
\]
Alternatively, $\vartheta\in T_{\Theta}^{\flat}(\theta_0)$ if and only if for every sequence $\{k_m\}$ diverging to $\infty$, there exists a sequence $\{\theta_m\}$ in $\Theta$ such that $k_m(\theta_m - \theta_0) \to \vartheta$.
Finally, if $\Theta$ is finite-dimensional, the Clarke tangent of $\Theta$ at $\theta_0$ is defined 
\[
C_{\Theta}(\theta_0) ~:=~ \liminf_{\Theta\ni\vartheta\to\theta_0}\, T_{\Theta}(\vartheta).
\]
(see Proposition 2.56 of~\cite{bonnans2013perturbation}.) It is not difficult to show that indeed these sets are cones. These cones are closed, the cone $C_{\Theta}(\theta_0)$ is always convex, and the inclusions $C_{\Theta}(\theta_0)\subseteq T_{\Theta}^{\flat}(\theta_0) \subseteq T_{\Theta}(\theta_0)$ hold for all $\Theta$ and $\theta_0\in\Theta$. If the cones $T_{\Theta}(\theta) \neq T_{\Theta}^{\flat}(\theta_0)$, then $\mathbb{X}_n(\cdot)$ may not converge in any reasonable sense, and hence, $r_n(\widehat{\theta}_n - \theta_0)$ may not converge in distribution.

Following~\citet[page 127]{aubin2009set}, we say $\Theta$ is {\em derivable} at $\theta_0\in\Theta$, if $T_{\Theta}^{\flat}(\theta_0) = T_{\Theta}(\theta_0)$, and we say $\Theta$ is {\em tangentially} regular at $\theta_0\in\Theta$, if $C_{\Theta}(\theta_0) = T_{\Theta}(\theta_0)$. \cite{geyer1994asymptotics} and~\cite{shapiro2000asymptotics} call such conditions {\em Chernoff regular} and {\em Clarke regular}, respectively. Note that a tangentially regular set $\Theta$ is also derivable. If $\Theta$ is derivable or Chernoff regular at $\theta_0$, then $\mathbb{X}(\cdot)$ in~\eqref{eq:limiting-distribution} can be taken to be
\[
\mathbb{X}(u) := \begin{cases}0, &\mbox{if }u\in T_{\Theta}(\theta_0),\\
+\infty, &\mbox{otherwise.} \end{cases}
\]
\cite{shapiro2000asymptotics} uses tangential regularity to study local $M$-estimators.
Below we provide a few properties and examples of tangent cones.
\begin{enumerate}
    \item For any $\Theta\subseteq\mathbb{R}^d$, if $\theta_0\in\mbox{int}(\Theta)$, i.e., $\theta_0$ is in the interior of $\Theta\subseteq\mathbb{R}^d$, then $T_{\Theta}(\theta_0) = \mathbb{R}^d$. This follows because there exists $\bar{\varepsilon} > 0$ such that $\theta_0 + \varepsilon u\in \Theta$ for every $\varepsilon\in[0, \bar{\varepsilon}]$ and $u\in S^{d-1}$ (i.e., $\|u\|_2 = 1$), which implies $S^{d-1}\subset T_{\Theta}(\theta_0)$ and because $T_{\Theta}(\theta_0)$ is a cone, this implies the result. 
    \item Suppose $g_1, \ldots, g_p$ are continuously differentiable real-valued functions on $\Omega$. If
    \[
    \Theta := \{\theta\in\Omega:\, g_j(\theta) \ge 0,\; 1 \le j\le p\},
    \]
    then, setting $\mathcal{I}(\theta_0) := \{j\in\{1, 2, \ldots, p\}:\, g_j(\theta_0) = 0\}$ as the subset of active constraints at $\theta_0\in\Theta$, we have
    \begin{equation}\label{eq:tangent-cone-derivation}
    T_{\Theta}(\theta_0) \subseteq \{\vartheta\in\Omega:\, \langle \nabla g_j(\theta_0), \vartheta\rangle \ge 0,\; j\in I(\theta_0)\},
    \end{equation}
    where $\nabla g_j(\cdot)$ is the derivative of $\theta\mapsto g_j(\theta)$. Moreover, if the {\em constraint qualification assumption} holds, i.e., 
    \[
    \exists v_0\in \Omega,\quad\mbox{such that}\quad \langle \nabla g_j(\theta_0), v_0\rangle > 0,\quad\mbox{for all}\quad j\in I(\theta_0),
    \]
    then equality holds in~\eqref{eq:tangent-cone-derivation}. See Eq. (2.19)--(2.20) of~\cite{shapiro2000asymptotics},~\citet[Proposition 4.3.7]{aubin2009set}, and~\citet[Lemma 12.2, Section 12.6]{Nocedal2006} for details. Also, see Proposition 2.61 of~\cite{bonnans2013perturbation} for the case where $g_j$'s are convex. As a special example of statistical importance, if 
    \[
    \Theta = \{\theta\in\mathbb{R}^d:\, \langle a_j, \theta\rangle \le b_j,\, 1\le j\le p,\, \langle a_j, \theta\rangle = b_j,\, p+1\le j\le q\},
    \]
    then, under Slater condition,
    \[
    T_{\Theta}(\theta_0) = \{\vartheta\in\mathbb{R}^d:\, \langle a_j, \vartheta\rangle \le 0,\, j\in \mathcal{I}(\theta_0),\, \langle a_j, \vartheta\rangle = 0,\, p+1 \le j\le q\}.
    \]
    \item Proposition 4.2.1 of~\cite{aubin2009set} states that all closed convex subsets $\Theta\subseteq\Omega$ with at least two elements are tangentially regular (and hence, also derivable) at every $\theta_0\in\Theta$. Also, see Proposition 2.57 of~\cite{bonnans2013perturbation}.
    \item (Example of derivable but not tangentially regular set.) Suppose $\Theta = \{(\theta_1, \theta_2)\in\mathbb{R}^2:\,|\theta_1| = |\theta_2|\}$. Then with $\theta_0 = (0, 0)$,
    \[
    T_{\Theta}(\theta_0) = T_{\Theta}^{\flat}(\theta_0) = \Theta,\quad\mbox{but}\quad C_{\Theta}(\theta_0) = \{\theta_0\}.
    \]
\end{enumerate}
In our discussion above, we considered the limit of $r_n(\widehat{\theta}_n - \theta_0)$ with a scalar $r_n$, which resulted in tangent cones in the limit. If $r_n$ is chosen to be a sequence of matrices, then more complications can be expected in finding the limiting set of constraints; see~\citet[Section 2]{pflug1995asymptotic}. 

In conclusion, identifying the limits of $\mathbb{D}_n(\cdot), \mathbb{S}_n(\cdot)$ and the tangent cones of $\Theta$ at $\theta_0$, yields the limiting distribution of $r_n(\widehat{\theta}_n - \theta_0)$. It is worth noting that the proof strategy above does not depend on the structure of $\mathbb{M}_n(\cdot)$ as an average. Moreover, the discussion of tangent cones highlights that the limiting distribution at least in the presence of constraints can significantly depend on $\theta_0$ (and hence, also on $P$). Especially, in such cases, traditional Wald, bootstrap, and subsampling approaches can fail to yield uniformly valid inference.
\subsection{Common Classes of Limiting Distributions of $M$-estimators}\label{sec:common-classes}
From~\eqref{eq:limiting-distribution}, it is clear that one can define $M$-estimators so that the limiting distribution statement~\eqref{eq:limiting-distribution} holds with a specified choice of arbitrary $\mathbb{D}(\cdot), \mathbb{S}(\cdot)$, as long as $\mathbb{D}(u) \ge 0$ for all $u\in T_{\Theta}(\theta_0)$ with $\mathbb{D}(0) = 0$, and $\mathbb{S}(\cdot)$ is a mean zero stochastic process. However, with some structure on the objective function of $M$-estimand (commonly encountered in statistical applications), one can find three most commonly occurring classes of limiting distributions. Pflug's classification is with respect to the structure/nature of the mean zero stochastic process $\mathbb{S}(\cdot)$. In all classes, the deterministic drift function $\mathbb{D}(\cdot)$ can be any non-negative function with a global minimum at zero and $\mathbb{X}(\cdot)$ is the indicator function of some cone.
\subsubsection{Class I: Linear Stochastic Processes}
The first class of limiting distributions has $\mathbb{S}(u) = \langle u, Y\rangle$ for some mean zero random vector $Y$. Most commonly, but not always, $Y\sim N(0, \Sigma)$. In general, it can be any infinitely divisible random vector. This class of limiting distributions arise when the loss function $\phi\mapsto L(\phi, W)$ is (almost everywhere) differentiable at $\theta_0$. Under differentiability, we get
\[
\mathbb{S}_n(u) \approx \frac{\rho_n}{r_n}\langle \nabla{\mathbb{M}}_n(\theta_0) - \nabla{\mathbb{M}}(\theta_0), \, u\rangle,
\]
where $\nabla{\mathbb{M}}_n(\cdot)$ and $\nabla{\mathbb{M}}(\cdot)$ are the gradients of $\mathbb{M}_n(\cdot)$ and $\mathbb{M}(\cdot)$, respectively. Note that for constrained problems, $\nabla \mathbb{M}(\theta_0)$ need not be zero. If $\rho_n/r_n$ is chosen so that $(\rho_n/r_n)(\nabla{\mathbb{M}}_n(\theta_0) - \nabla{\mathbb{M}}(\theta_0)))$ converges in distribution to $Y$, then $\mathbb{S}(u) = \langle Y, u\rangle$. In case, $\mathbb{M}_n(\cdot)$ is an average, $Y$ can be any infinitely divisible random vector.

The following are some examples where linear stochastic processes occur in the limiting distribution.
\begin{example}[Constrained Mean Estimation] \label{ex:mean}
Suppose $W_1, W_2, \hdots, W_n$ are IID random vectors from some distribution with unknown mean vector $\theta_0$ and unknown covariance matrix $\Sigma$. Suppose that $\theta_0$ is known to lie in $\Theta \subset \R^p$, assumed to be closed and non-empty. The natural estimator $\widehat{\theta}_n$ of the mean vector is given by
\begin{equation*}
    \widehat{\theta}_n = \argmin_{\phi \in \Theta} \frac{1}{n}\sum_{i=1}^n \norm{W_i - \phi}_2^2.
\end{equation*}
If $\Theta$ is derivable or Chernoff regular at $\theta_0$, then it follows from~\cite{geyer1994asymptotics} and~\cite{shapiro2000asymptotics} that
\begin{equation*}
    n^{1/2}(\widehat{\theta}_n - \theta_0) ~\overset{d}{\to}~ \argmin_{u \in \R^p} \: -2\langle u, Y\rangle + \mathbb{D}(u) + \mathbb{X}(u),
\end{equation*}
if $n^{-1/2}\sum_{i=1}^n (W_i - \theta_0) \overset{d}{\to} Y$, where $\mathbb{D}(u) = \|u\|_2^2$, and $\mathbb{X}(u)$ is the indicator function of the tangent cone $T_{\Theta}(\theta_0)$. If $W_i$'s satisfy the Lindeberg condition, then $Y$ would be a mean zero normal random vector. For a specific example, if $\Theta = [0, \infty) \subseteq\mathbb{R}$, then $T_{\Theta}(\theta_0) = [0, \infty)$ if $\theta_0 = 0$ and $T_{\Theta}(\theta_0) = \mathbb{R}$ if $\theta_0 \neq 0$. Hence, $n^{1/2}(\widehat{\theta}_n - \theta_0) \overset{d}{\to} \max\{N(0, \sigma^2), 0\}$ if $\theta_0 = 0$ and $n^{1/2}(\widehat{\theta}_n - \theta_0) \overset{d}{\to} N(0, \sigma^2)$, otherwise. This example can be generalized to the case where $\widehat{\theta}_n$ is the $M$-estimator with arbitrary loss function $L(\phi, W)$ satisfying $\mathbb{E}[\nabla L(\theta_0, W)] = 0$; see Assumption A and the results of~\cite{geyer1994asymptotics}.
\end{example}
\begin{example}[Median/LAD Regression] \label{ex:L1}
    Suppose $(x_i, Y_i), 1\le i\le n$ are $\mathbb{R}^{p+1}$ dimensional random vectors with non-stochastic covariates $x_1, \ldots, x_n$, satisfying the linear model 
    \begin{equation*}
             Y_i = \theta^{(0)}_0 + \theta^{(1)}_0 x_{1i} + \cdots + \theta^{(p)}_0 x_{pi} + \xi_i,
    \end{equation*}
    where $\theta_0 = (\theta_0^{(0)}, \theta_0^{(1)}, \hdots, \theta_0^{(p)})^\top$ is a vector of unknown parameters and $\{\xi_i\}$ are unobservable IID errors with median 0. The median/least absolute deviation/$L_1$ regression estimator of $\theta_0$ is defined as 
    \begin{equation*}
        \widehat{\theta}_n = \argmin_{\phi\in\mathbb{R}^{p+1}} \sum_{i=1}^{n} | Y_i - \langle \phi, x_i\rangle|,
    \end{equation*}
    where $\langle \phi, x_i\rangle = \phi^{(0)} + \sum_{j=1}^p \phi^{(j)}x_{ji}.$
    Set $F_i(e) = \mathbb{P}(\xi_i \le e)$ and for some sequence of constants $a_n$ and $t > 0$, define
    \[
    \Psi_{ni}(t) = \int_0^t n^{1/2}(F_i(s/a_n) - F_i(0))ds,\quad\mbox{and}\quad \Psi_{ni}(-t) = \int_{-t}^{0} n^{1/2}(F_i(s/a_n) - F_i(0))ds
    \]
    \cite{knight1999asymptotics} obtained the limiting distribution of $\widehat{\theta}_n$ under the following assumptions, which allow for zero/non-existent density of $F_i(\cdot)$: 
    \begin{enumerate}[label=(LAD-\arabic*)]
    \item $\xi_i$'s are independent random variables with a continuous distribution function and $F_i(0) = 1/2$;\label{continuity-LAD}
    \item for some positive definite matrix $C$, $n^{-1}\sum_{i=1}^n \langle \phi, x_i\rangle^2 \to \langle \phi, C\phi\rangle$ as $n\to\infty$ for all $\phi\in\mathbb{R}^{p+1}$, and $n^{-1}\max_{1\le i\le n}\|x_i\|_2 \to 0$ as $n\to\infty$; and\label{design-LAD} \item there exists a non-negative strictly convex function $\mathbb{D}(\cdot)$ such that $n^{-1}\sum_{i=1}^n \Psi_{ni}(\langle u, x_i\rangle) \to \mathbb{D}(u)$ as $n\to\infty$ for all $u\in\mathbb{R}^{p+1}$. \label{differentiability-LAD}
    \end{enumerate}
    The generality of these assumptions is discussed at length in~\cite{knight1999asymptotics,knight-1998}. Under assumptions~\ref{continuity-LAD}--\ref{differentiability-LAD}, Theorem 1 of~\cite{knight1999asymptotics} states
    \begin{equation*}
        a_n(\widehat{\theta}_n - \theta_0) ~\overset{d}{\to}~ \argmin_{u \in \R^{p+1}} \; -\langle u, Y\rangle + \mathbb{D}(u),
    \end{equation*}
    where $Y\sim N(0, C) \in \mathbb{R}^{p+1}$.
\end{example} 

\begin{example}[Bridge Estimators] \label{ex:L2}
Consider the linear regression model 
    \begin{equation*}
        Y_i = \theta_0^{(1)} x_{1i} + \cdots + \theta_0^{(p)} x_{pi} + \xi_i,
    \end{equation*}
    where $\theta_0 = (\theta_0^{(1)}, \hdots, \theta_0^{(p)})^\top$ is a vector of unknown parameters, covariates are non-stochastic, and $\{\xi_i\}$ are unobservable IID errors with mean 0 and variance $\sigma^2$. (Following~\cite{fu2000asymptotics}, we assume zero intercept.) The parameters are estimated by minimizing a penalized least squares criterion, i.e., we have 
    \begin{equation*}
        \widehat{\theta}_n = \argmin_{\phi\in\mathbb{R}^p} \sum_{i=1}^{n} (Y_i - \phi^{(1)} x_{1i} - \cdots - \phi^{(p)} x_{pi})^2 + \lambda_n \sum_{j=1}^{p} |\phi^{(j)}|^{\mu}
    \end{equation*}
    for a given $\{ \lambda_n \}$ and $\mu > 0$. The ridge and LASSO regression estimators are special cases when $\mu = 2$ and $\mu = 1$ respectively. Under assumption~\ref{design-LAD} with a non-negative definite matrix $C$, if $\lambda_n/n^{\min\{\mu, 1\}/2} \to \lambda_0\in(0,\infty)$ as $n\to\infty$, then Theorems 2 and 3 of~\cite{fu2000asymptotics} state that
    \begin{equation*}
        n^{1/2}(\widehat{\theta}_n - \theta_0) ~\overset{d}{\to}~ \argmin_{u \in \R^p} \; -2\langle u, Y\rangle + \mathbb{D}(u),
    \end{equation*}
    where $Y \sim N(0,\sigma^2 C)$ and $\mathbb{D}(u) = \langle u, Cu\rangle + \lambda_0 \sum_{j=1}^p h_{\mu}(u^{(j)}, \theta_0^{(j)}),$ with $h_{\mu}:\mathbb{R}^2\to\mathbb{R}$ defined as
    \begin{equation*}
        h_{\mu}(a, b) = 
        \begin{cases}
            a|b|^{\mu - 1}\sign{(b)} & \text{if } \mu > 1, 
                \\
             a \,\sign{(b)}\mathbf{1}\{b \neq 0\} + |a|\mathbf{1}\{b = 0\}  & \text{if } \mu = 1, \\
            |a|^{\mu} \, \mathbf{1}\{b = 0\}  & \text{if } \mu < 1.
        \end{cases}
    \end{equation*}
    Unlike the previous examples, the deterministic function $\mathbb{D}(\cdot)$ can be convex or non-convex depending on what $\theta_0$ is.
\end{example}

\subsubsection{Class II: Gaussian Process} 
The second class of limiting distributions has $\mathbb{S}(\cdot)$ as a (non-linear) Gaussian process. This occurs in the case where the loss function $L(\phi, W)$ is not ``smooth'' in $\phi$; a prototypical example in~\cite{pflug1995asymptotic} is $L(\phi, W) = H(W)\mathbf{1}\{V(\phi, W) \ge 0\}$ for some smooth function $\phi\mapsto V(\phi, W)$. Many of the cube-root estimators of~\cite{kim1990cube} have Gaussian processes in the limiting distribution. For this scenario, the mean zero stochastic process is $\mathbb{S}(u) = \mathrm{GP}_{\Sigma}(u)$, indexed by the covariance operator $\Sigma(\cdot, \cdot)$. Note that the classical two-sided Brownian motion (or Wiener process) is a special case. \cite{pflug1995asymptotic} considers Gaussian processes with covariance operators of a specific form 
\begin{equation}\label{eq:covariance-Pflug-form}
\Sigma(u, v) = \gamma\int \frac{y^2}{2}(|\langle s, u\rangle| + |\langle s, v\rangle| - |\langle s, u-v\rangle|)d\nu(y, s),
\end{equation}
for some $\gamma > 0$ and some probability measure $\nu$ on $\mathbb{R}\times S^{p-1}$. (Note that $\min\{|u|, |v|\}\mathbf{1}\{uv \ge 0\} = (1/2)(|u| + |v| - |u-v|)$.)
For simplicity and easy connection to statistical applications, we consider general (non-linear) Gaussian processes.

The following are some examples where (non-linear) Gaussian processes occur in the limiting distribution. 
\begin{example}[Shorth estimator] \label{ex:shorth}
    Suppose $W_1, W_2, \hdots, W_n$ are independent observations sampled from a distribution $P_0$ on the real line, with a Lebesgue density $p_0$. A popular robust location estimator in this context is the traditional shorth estimator $\widehat{\theta}_n$ \citep{pensia2019mean}, which is defined as the center of the shortest interval containing at least $n/2$ observations. More precisely,  
\begin{align*}
    \widehat{r}_n ~&:=~ \inf \Big\{ r : \sup_{x} \sum_{i=1}^n \mathbf{1}\{W_i\in [x - r, \, x + r ]\} \geq n/2 \Big\}, \\
    \widehat{\theta}_n ~&:=~ \argmax_{x} \: \sum_{i=1}^n \mathbf{1}\{W_i\in [x - \widehat{r}_n, \, x + \widehat{r}_n]\}.
\end{align*} 
If $[\theta_0 - r_0, \, \theta_0 + r_0]$ is the unique shortest interval with $P_0$ measure at least half, then assuming that $p_0(\theta_0 + \eta r_0) > 0$ for $\eta\in\{\pm 1\}$, and that $c_0 := p'_0(\theta_0 - r_0) - p'_0(\theta_0 +r_0) > 0$, Example 6.1 of~\cite{kim1990cube} states that  
\begin{equation*}
    n^{1/3}(\widehat{\theta}_n - \theta) ~\overset{d}{\to}~ \argmin_{u \in \R} \: \mathrm{GP}_{\Sigma}(u) + \mathbb{D}(u),
\end{equation*}
where $\Sigma(u,v) = 2p_0(\theta_0 + r_0)\min\{|u|, |v|\}\mathbf{1}\{uv > 0\}$, and $\mathbb{D}(u) = c_0u^2/2$.
\end{example}

\begin{example}[Least Median of Squares] \label{ex:median}

Suppose $(X_i, Y_i)$, $1 \leq i \leq n$ are $\R^{p+1}$ dimensional random vectors with stochastic covariates $X_1, \hdots, X_n$ sampled independently from a distribution $P_0$ on $\R^p$, satisfying the linear model $Y_i = \langle \theta_0, X_i \rangle + \xi_i,$
where $\theta_0\in\mathbb{R}^p$ is an unknown parameter and $\{\xi_i\}$ are unobservable errors. \cite{rousseeuw1984least} defined the least median of squares estimator of $\theta_0$ as
\begin{equation*}
    \widehat{\theta}_n = \argmin_{\phi\in\mathbb{R}^p} \, \median\limits_{1\leq i \leq n} (Y_i - \langle \phi, X_i \rangle)^2.
\end{equation*} 
\cite{kim1990cube} obtained the limiting distribution of $\widehat{\theta}_n$, under the following assumptions:

\begin{enumerate}[label=(LMS-\arabic*)]
\item $X_i$ and $\xi_i$ are independently distributed. \label{indep.-LMS}
\item $X_i$ has finite second moment and $Q := \E_{P_0}[X_i X_i^\top]$ is positive definite.
\item $\xi_i$ has a bounded, symmetric density $p_0$ that decreases away from zero, and $p_0'(m_0) < 0$, where $m_0$ is the unique median of $|\xi_i|$. \label{density-LMS}
\end{enumerate}
Under assumptions~\ref{indep.-LMS}-\ref{density-LMS}, Example 6.4 of \cite{kim1990cube} shows that  
\begin{equation*}
    n^{1/3}(\widehat{\theta}_n - \theta_0) ~\overset{d}{\to}~ \argmin_{u \in \R^p} \: \mathrm{GP}_{\Sigma}(u) + \mathbb{D}(u),
\end{equation*}
where $\Sigma(u, v) = p_0(m_0) \: \mathbb{E}_{P_0} \left[|\langle X, u\rangle| + |\langle X, v\rangle| - |\langle X, u-v\rangle| \right]$ and $\mathbb{D}(u) =  -p_0'(m_0)\, \langle u, Qu \rangle$.
\end{example}

\begin{example}[Mode Estimation] \label{ex:mode}
Suppose $W_1, W_2, \hdots, W_n$ are IID random variables drawn from a unimodal density $p_0$ on $\mathbb{R}$. Assume that $p_0$ is continuous and achieves its maximum at the unique point $\theta_0$. \citet{venter1967estimation} proposed the mode estimator $\widehat{\theta}_n$ defined by
\begin{gather*}
    K_n = \argmin_{r_n + 1 \leq j \leq n-r_n} W_{j+r_n} - W_{j-r_n}, \qquad \quad
    \widehat{\theta}_n = W_{K_n},
\end{gather*}
where $\{r_n\}$ is an appropriately chosen sequence of integers. This estimator can be considered the maximizer of spacings (or nearest neighbor) based density estimator. If $p_0$ is thrice continuously differentiable at $\theta_0$ with $\kappa_1 = p_0'(\theta_0) = 0$, $\kappa_2 = -p_0''(\theta_0) > 0$, and $\kappa_3 = p_0'''(\theta_0)$, so that 
\begin{equation*}
    p_0(x) = p_0(\theta_0) - \frac{1}{2} \kappa_2 (x-\theta_0)^2 + \frac{1}{6} \kappa_3 (x-\theta_0)^3 + o(|x-\theta_0|^3) \; \; \text{as } \; x \rightarrow \theta_0,
\end{equation*} 
then, for $r_n \sim An^{\mu}$ with $4/5 \leq \mu \leq 7/8$, Theorem 3b of \citet{venter1967estimation} states that 
\begin{equation*}
    2^{-1/3} (A\kappa_2)^{2/3} p_0(\theta_0)^{-1} n^{(2 \mu-1)/3} (\widehat{\theta}_n - \theta_0) ~\overset{d}{\to}~  \argmin_{u \in \R} \: \mathbb{S} (u) + \mathbb{D}(u),
\end{equation*}
where $\mathbb{S}(u) = \mathrm{GP}_{\Sigma}(u)$ with $\Sigma(u, v) = \min\{|u|, |v|\}\mathbf{1}\{uv > 0\}$, and $\mathbb{D}(u) = u^2 - c_0u\mathbf{1}\{\mu = 7/8\}$ for $c_0 = 3^{-1} 2^{-1/3} A^{8/3} \kappa_2^{-1/3} \kappa_3 p_0(\theta_0)^{-4}$.
\end{example}

\subsubsection{Class III: Generalized Poisson Hyperplane Process} 
The third and final class of limiting distributions we consider has $\mathbb{S}(\cdot)$ as a centered Poisson hyperplane process. As explained in~\cite{pflug1995asymptotic}, such a process occurs in non-smooth $M$-estimation problems with constraints. The (generalized) Poisson hyperplane process denoted by $\text{Poisson}(\gamma, {\nu}) (u)$, indexed by $\gamma \ge 0$ and a probability measure $\nu$ on $\mathbb{R}^{1+p}$, is formally defined as follows.
Let $\{\tau_i\}_{i \geq 1}$ be the jump times of a Poisson process with parameter $\gamma$, i.e., $\tau_1$ and $\tau_{i+1} - \tau_{i}, i \geq 1$ are IID exponential random variables with expectation $1/\gamma$. If $(V_i, U_i)_{i \geq 1}$ is an independent sequence of random vectors distributed according to $\nu$, then the Poisson hyperplane process is given by
\begin{equation*}
    \mbox{Poisson}(\gamma, \nu)(u) = \sum_{i=1}^{\infty} V_i \mathbf{1}\{ \tau_i \leq \langle u, U_i\rangle\}, \qquad u \in \R^p.
\end{equation*}

Then the mean zero stochastic process occurring in the limiting distribution is the centered Poisson hyperplane process $$\mathbb{S}(u) ~=~ \mathrm{cPoisson}(\gamma, \nu)(u) ~:=~ \mbox{Poisson}(\gamma, \nu)(u) - \mathbb{E}[\mbox{Poisson}(\gamma, \nu)(u)].$$ Note that $\mathbb{E}[\mbox{Poisson}(\gamma, \nu)(u)] = \gamma\int y(\langle u, s\rangle)_+ d\nu(y, s)$.
\begin{remark}
The traditional Poisson hyperplane process, as described by \citet{pflug1995asymptotic}, restricts the measure ${\nu}$ to be supported only on $\R \times S^{p-1}$. We generalize the support to include additional limiting processes that arise in this context.
\end{remark}

\begin{example}[Linear Regression: Discontinuous densities] \label{ex:regression_discont}

Suppose $(Y_i, X_i)_{i \geq 1}$ is an IID sequence of vectors in $\R \times \R^p$ satisfying the linear model
\begin{equation*}
    Y_i = \langle X_i, \theta \rangle + \xi_i,
\end{equation*}
parametrized by $\theta \in \Theta$, a compact, convex subset of $\R^p$, $\{X_i\}_{i \geq 1}$ are stochastic covariates with common c.d.f. $F_X$, compact support $\mathcal{X}$ that is independent of $\theta$ and the error $\xi_i$ has conditional density $f(\xi| X_i, \theta)$, which has a jump at zero:
\begin{align*}
    &\lim_{\xi \uparrow 0} f(\xi| x, \theta) = q(x, \theta), \quad \lim_{\xi \downarrow 0} f(\xi| x, \theta) = p(x, \theta), \\
    &p(x, \theta) - \delta > q(x, \theta) > c, \quad \delta,\, c  > 0, \; \forall \: x \in \mathcal{X}, \; \forall \: \theta \in \Theta.
\end{align*}
Denote by $\theta_0$ the true value of the parameter and assume $\theta_0 \in \emph{int}(\Theta)$. The maximum likelihood estimator of $\theta_0$ is defined as
\begin{equation*}
    \widehat{\theta}_n ~:=~ \argmax_{\theta\in\Theta} \,\prod_{i=1}^n f(Y_i - \langle X_i, \theta\rangle | X_i, \theta) dF_X(X_i).
\end{equation*} 
 Under mild regularity assumptions on the conditional density, Theorem 3 of \citet{chernozhukov2002likelihood} states that
\begin{equation*}
    n(\widehat{\theta}_n - \theta_0) ~\overset{d}{\to}~ \argmin_{u \in \R^p} \: \mathrm{cPoisson} \left(2, (\nu_1 + \nu_2)/2 \right) (u) + \mathbb{D}(u),
\end{equation*} 
provided that the limiting process attains a unique minimum almost surely. 
Here $\nu_1$ and $\nu_2$ are the probability measures induced on $\R^{1+p}$ by the random vectors $(\log \: p(X)/q(X), \: p(X) X)$ and $(-\log \: p(X)/q(X), -q(X)X )$ respectively, with $p(X) = p(X, \theta_0)$ and $q(X) = q(X, \theta_0)$. The deterministic drift function is 
\begin{equation*}
    \mathbb{D}(u) = \E \left[ \left\langle u, X \right \rangle \Big( q(X) - p(X) \Big) + \left(\left\langle u, X \right \rangle \right)_+
    \left( \log \frac{p(X)}{q(X)} \big(p(X) + q(X)\big) \right) \right].
\end{equation*}

\end{example}

\section{Sufficient Conditions for Symmetry}\label{sec:sufficiency}
In this section, we first give sufficient conditions for symmetry of $\argmin_{u} \mathbb{Z}(u)$ in the general case without specific structure on the stochastic process. We then discuss implications for the three special classes described in Section~\ref{sec:common-classes} by providing equivalent criteria for each one and revisiting the accompanying examples. 
\begin{lemma}\label{lem:sufficiency-general}
    Suppose $\Omega$ is a symmetric set (i.e., $u\in\Omega$ implies $-u\in\Omega$) and $u\mapsto \mathbb{Z}(u), u\in \Omega$ is a stochastic process taking values in $\overline{\mathbb{R}} = [-\infty, \infty]$. Suppose the stochastic processes $u\mapsto \mathbb{Z}(u)$ and $u\mapsto \mathbb{Z}(-u)$ have the same law.\footnote{In the case $\Omega = \mathbb{R}^p, p\ge1$, this is equivalent to the assertion that for any $k\ge1$ and $u_1, \ldots, u_k\in\mathbb{R}^p$, $(\mathbb{Z}(u_1), \ldots, \mathbb{Z}(u_k))\overset{d}{=} (\mathbb{Z}(-u_1), \ldots, \mathbb{Z}(-u_k))$.} If $\mathbb{Z}(\cdot)$ is almost surely uniquely minimized at a random variable $W\in\Omega$, then $W$ is symmetric around zero in the sense that $W \overset{d}{=} -W$. 
\end{lemma}
\begin{proof}
    Note that for any Borel measurable set $A$,
    \begin{align*}
    \mathbb{P}(W \in A) &= \mathbb{P}\left(\inf_{u\in A}\, \mathbb{Z}(u) < \inf_{u\notin A}\, \mathbb{Z}(u)\right),\\
    \mathbb{P}(W \in -A) &= \mathbb{P}\left(\inf_{u\in -A}\,\mathbb{Z}(u) < \inf_{u\notin -A}\, \mathbb{Z}(u)\right)\\
    &= \mathbb{P}\left(\inf_{u\in A}\,\mathbb{Z}(-u) < \inf_{u\notin A}\, \mathbb{Z}(-u)\right).
    \end{align*}
    Because $u\mapsto \mathbb{Z}(u)$ and $u\mapsto \mathbb{Z}(-u)$ have the same law, we conclude that
    \[
    \left(\inf_{u\in A} \mathbb{Z}(u),\,\inf_{u\notin A}\mathbb{Z}(u)\right) \overset{d}{=} \left(\inf_{u\in A}\mathbb{Z}(-u),\,\inf_{u\notin A} \mathbb{Z}(-u)\right).
    \]
    This, in turn, implies $\mathbb{P}(W \in A) = \mathbb{P}(-W \in A)$ for all Borel measurable sets $A$, i.e., $W \overset{d}{=} -W$.
\end{proof}

For notational convenience, we call a stochastic process $u\mapsto \mathbb{Z}(u)$ {\em even} if $u\mapsto \mathbb{Z}(u)$ and $u\mapsto \mathbb{Z}(-u)$ share the same law. (Note that if $\mathbb{Z}(\cdot)$ is non-stochastic, this definition reduces to that of an even function.)

In the cases of interest to us where $\mathbb{Z}(u)$ decomposes as $\mathbb{D}(u) + \mathbb{S}(u) + \mathbb{X}(u)$, the sufficient condition of Lemma~\ref{lem:sufficiency-general} is equivalent to requiring even-ness of the deterministic function $\mathbb{D}(u)$, symmetry of the mean zero stochastic process $\mathbb{S}(u),$ and the symmetry of the set to which $\mathbb{X}(\cdot)$ is the indicator function. Given the importance of this result in the context of $M$-estimation, we state this as a theorem.
\begin{theorem} \label{thm:suff}
Suppose $\Theta^\star$ is a symmetric set (i.e., $u\in\Theta^\star$ implies $-u\in\Theta^\star$), and a non-stochastic map $\mathbb{D}(\cdot)$ defined on $\Theta^{\star}$ is even (i.e., $\mathbb{D}(-u) = \mathbb{D}(u)$ for all $u\in\Theta^{\star}$). Additionally, if $u\mapsto \mathbb{S}(u)$ is a mean zero {\em even} stochastic process, then $W = \argmin_{u\in\Omega} \mathbb{D}(u) + \mathbb{S}(u) + \mathbb{X}(u)$ is a symmetric (around zero) random variable. (Here $\mathbb{X}(u) = 0$ for all $u\in\Theta^\star$ and $\mathbb{X}(u) = +\infty$ for all $u\notin\Theta^\star$.)
\end{theorem}

In the following, we provide sufficient conditions for the even-ness of the mean zero stochastic process $\mathbb{S}(\cdot)$. We will throughout assume the following. 
\begin{enumerate}[label=\bf(A\arabic*)]
    \item The deterministic function $u\mapsto \mathbb{D}(u) + \mathbb{X}(u)$ is even.\label{assump:even-ness-sufficiency} 
\end{enumerate}
Observe that if $\mathbb{X}(\cdot)$ is the indicator function of a set $\Theta^\star$ and $\mathbb{D}(\cdot)$ is real-valued, then assumption~\ref{assump:even-ness-sufficiency} is equivalent to the symmetry of $\Theta^\star$ and even-ness of $\mathbb{D}(\cdot)$ on $\Theta^\star$. In the absence of constraints, assumption~\ref{assump:even-ness-sufficiency} is not difficult to prove or disprove in the examples considered. 
\setcounter{example}{0}
\setcounter{section}{2}
\begin{example}[Constrained Mean Estimation (Revisited)]
    In constrained mean estimation problem, $\mathbb{D}(u) = \|u\|_2^2$ which is an even function. Hence, assumption~\ref{assump:even-ness-sufficiency} holds if and only if the tangent cone is a symmetric set. Note that if $\theta_0$ belongs to the interior of the constraint set $\Theta$, then the tangent cone is trivially symmetric. From~\eqref{eq:tangent-cone-derivation}, it follows that for constraint sets defined via (non-linear) inequality and equality constraints need not be symmetric if $\theta_0$ is a point on the boundary of $\Theta$. In fact, for linear inequality and equality constraints, symmetry of the tangent cone holds if and only if the tangent cone is a subspace.
\end{example}
\begin{example}[Median/LAD Regression (Revisited)]
    In the LAD regression problem, $\mathbb{X}(u)\equiv 0$ and $\mathbb{D}(u)$ is the limit of $n^{-1}\sum_{i=1}^n \Psi_{ni}(\langle u, x_i\rangle)$. A simple sufficient condition for symmetry of $\mathbb{D}(u)$ is
    \[
    \lim_{s\downarrow 0}\max_{1\le i\le n}\frac{F_i(s)-F_i(0)}{F_i(0) - F_i(-s)} ~=~ 1.
    \]
    As a special case, if $F_i(s) - F_i(0) = \lambda_i\sign{(s)}|s|^{\gamma}L(|s|)$ for some $\gamma \in (0,\infty)$ and a slowly varying function $L(\cdot)$, then $\mathbb{D}(u)$ is the limit of $(\gamma + 1)^{-1}n^{-1}\sum_{i=1}^n \lambda_i|\langle u, x_i\rangle|^{\gamma + 1}$ (as $n\to\infty$), which is an even function of $u$. Note that $\gamma = 1$ is the standard/regular case where the errors have non-zero Lebesgue densities at zero. 
\end{example}
\begin{example}[Bridge Estimators (Revisited)]
    In case of bridge estimators, $\mathbb{X}(u) = 0$ and $\mathbb{D}(u) = \langle u, Cu\rangle + \lambda_0\sum_{j=1}^p h_{\mu}(u^{(j)}, \theta_0^{(j)})$, where
    \begin{equation*}
        h_{\mu}(a, b) = 
        \begin{cases}
            a|b|^{\mu - 1}\sign{(b)} & \text{if } \mu > 1, 
                \\
             a \,\sign{(b)}\mathbf{1}\{b \neq 0\} + |a|\mathbf{1}\{b = 0\}  & \text{if } \mu = 1, \\
            |a|^{\mu} \, \mathbf{1}\{b = 0\}  & \text{if } \mu < 1.
        \end{cases}
    \end{equation*}
    Recall that the penalty of the bridge estimator is non-convex for $\mu \in (0,1)$ and convex for $\mu \ge 1.$ Validity of assumption~\ref{assump:even-ness-sufficiency} shows an interesting dichotomy between these cases. The quadratic part of $\mathbb{D}(u)$ is always even. For $\mu < 1$, $a\mapsto h_{\mu}(a, b)$ is an even function of $a$, no matter the value of $b$. For $\mu \ge 1$, $a\mapsto h_{\mu}(a, b)$ is even if and only if $b = 0$. Hence, assumption~\ref{assump:even-ness-sufficiency} holds for all $\mu \in (0, 1)$, for all $\theta_0\in\mathbb{R}^p$ and for all $\mu\in(1, \infty)$ if and only if $\theta_0 = 0$. In other words, assumption~\ref{assump:even-ness-sufficiency} does not hold uniformly over the parameter space.  
\end{example}
\begin{example}[Shorth Estimator (Revisited)]
In case of the Shorth estimator, $\mathbb{D}(u) = c_0u^2/2$ and hence, assumption~\ref{assump:even-ness-sufficiency} holds. Note that if $c_0 = 0$, then the asymptotic distribution changes which could break assumption~\ref{assump:even-ness-sufficiency}.
\end{example}
\begin{example}[Least Median of Squares (Revisited)]
In case of the least median of squares estimator, $\mathbb{D}(u) = -p_0'(m_0)\langle u, Qu\rangle$ for a positive definite matrix $Q$. Hence, assumption~\ref{assump:even-ness-sufficiency} holds true.
\end{example}
\begin{example}[Mode Estimation (Revisited)]
In case of the mode estimator, $\mathbb{D}(u) = u^2 - c_0u\mathbf{1}\{\mu = 7/8\}$ for a tuning parameter $\mu$ of the estimation procedure. For $\mu\in[4/5, 7/8)$, the resulting estimator satisfies assumption~\ref{assump:even-ness-sufficiency}, but for $\mu = 7/8$, assumption~\ref{assump:even-ness-sufficiency} does not hold.
\end{example}
\begin{example}[Linear Regression: Discontinuous densities (Revisited)]
In case of linear regression with discontinuous error density, the deterministic drift $\mathbb{D}(u)$ is of the form $\mathbb{D}(u) = \mathbb{E}[\langle u, X\rangle R_1] + \mathbb{E}[(\langle u, X\rangle)_+R_2]$, for a random vector $X$ and two real-valued random variables $R_1, R_2$. Using the fact that $(x)_+ = (|x| + x)/2$ for all $x\in\mathbb{R}$, this can be equivalently rewritten as
\[
\mathbb{D}(u) = \mathbb{E}[|\langle u, X\rangle|R_2/2] + \mathbb{E}[\langle u, X\rangle(R_1 + R_2/2)].
\]
This implies that
\[
\mathbb{D}(u) - \mathbb{D}(-u) = 2\langle u, \mathbb{E}[X(R_1 + R_2/2)]\rangle,
\]
which is zero for all $u\in\mathbb{R}^p$ if and only if $\mathbb{E}[X(R_1 + R_2/2)] = 0.$
Hence, even-ness of $\mathbb{D}(\cdot)$ (or equivalently,~\ref{assump:even-ness-sufficiency}) is not guaranteed, in general.
\end{example}
\setcounter{section}{3}
\begin{corollary}[Class I: Linear Stochastic Process]\label{suff:NS}
    The stochastic process $\mathbb{S}(u) = \langle u, Y\rangle$ defined on $\mathbb{R}^p$ is {\em even} if and only if $Y \overset{d}{=} -Y$. Hence, if $Y$ is symmetric, class I limiting distributions are symmetric under~\ref{assump:even-ness-sufficiency}.
\end{corollary} 
In Examples~\ref{ex:mean}--\ref{ex:L2}, $\mathbb{S}(u) = \langle u, Y\rangle$ for a mean zero Gaussian random vector $Y$. Because all mean zero Gaussian random vectors are symmetric around zero, the limiting distributions are symmetric under~\ref{assump:even-ness-sufficiency}. It should be stressed here that because in Class I, $\mathbb{S}(u)$ arises as the limiting distribution of $\langle \nabla{\mathbb{M}}_n(\theta_0) - \nabla{\mathbb{M}}(\theta_0), u\rangle$, $Y$ can be any element of the class of infinitely divisible distributions, which means that examples exist where Class I distributions can be asymmetric. For a simple one, consider the centered Poisson limiting distribution of centered Binomial($n, 1/n$), which is asymmetric.

\begin{corollary}[Class II: Gaussian Process]\label{suff:GP}
    The (non-linear) Gaussian process $\mathbb{S}(u) = \emph{GP}_{\Sigma}(u)$ defined on $\mathbb{R}^p$ is {\em even} if and only if $\Sigma(u, v) = \Sigma(-u, -v)$  for all $u, v \in \R^p$. Hence, if $\Sigma(u, v) = \Sigma(-u, -v)$  for all $u, v$, then under~\ref{assump:even-ness-sufficiency}, class II limiting distributions are symmetric.    
\end{corollary} 
In Examples~\ref{ex:shorth}--\ref{ex:mode}, the mean zero Gaussian process appearing in the limiting distribution is even because they all share the form~\eqref{eq:covariance-Pflug-form} from Pflug's work. Because the function $R_s(u, v) = |\langle s, u\rangle| + |\langle s, v\rangle| - |\langle s, u - v\rangle|$ satisfies $R_s(u, v) = R_s(-v, -u)$ for all $u, v$, every covariance matrix of the form~\eqref{eq:covariance-Pflug-form} satisfies the assumptions of Corollary~\ref{suff:GP}, no matter what the measure $\nu(\cdot)$ is.

\begin{corollary}[Class III: Generalized Poisson Hyperplane Process]\label{suff:GPHP}
    The (generalized) Poisson hyperplane process $\emph{\mbox{Poisson}}(\gamma, \nu)(u)$ defined on $\mathbb{R}^p$ is {\em even} if $\nu(A, B) = \nu (A, -B)$ for all Borel measurable sets $A\subseteq \R, B \subseteq \R^p$. Hence, if $\nu(A, B) = \nu (A, -B)$ for all Borel measurable $A, B$, then the centered Poisson hyperplane process $\mathbb{S}(u) = \emph{\mbox{cPoisson}}(\gamma, \nu)(u)$ is {\em even} and under~\ref{assump:even-ness-sufficiency}, class III limiting distributions are symmetric.    
\end{corollary} 
Unlike the case for Classes I and II, for Class III, we do not know of a necessary condition for the even-ness of a  Poisson hyperplane process. Moreover, Example~\ref{ex:regression_discont} does {\em not} satisfy the assumptions of Corollary~\ref{suff:GPHP} and hence, we cannot confirm whether the even-ness holds true or not. 

\section{Necessary Conditions for Symmetry of Limiting Distribution: Linear Stochastic Processes} \label{Random Shift}\label{sec:necessity}
In Section~\ref{sec:sufficiency}, we provided simple sufficient conditions for symmetry of the limiting distributions of $M$-estimators. For most cases (e.g., classes I and II), assumption~\ref{assump:even-ness-sufficiency} is the bottleneck to verify the symmetry of the limiting distribution. If we can verify the necessity of assumption~\ref{assump:even-ness-sufficiency} for symmetry, then for classes I and II we have obtained necessary and sufficient conditions for symmetry of the limiting distributions. In what follows, we show the necessity of assumption~\ref{assump:even-ness-sufficiency} for class I limiting distributions. The technicalities included even for this simple case make us believe that proving the necessity of~\ref{assump:even-ness-sufficiency} in general is a much harder problem. One other assumption of Theorem~\ref{thm:suff} is that the minimizer of the stochastic process is almost surely unique. In Section~\ref{sec:exist-arg-min}, we provide sufficient conditions for the existence and uniqueness of the minimizer of class I stochastic processes. Although it is the most common, we do not assume $Y$ to be a mean zero Gaussian.

Consider the stochastic process 
\begin{equation}\label{eq:definition-Z}
    \mathbb{Z}(u) = \mathbb{D}(u) + \mathbb{X}(u) +\langle u, Y\rangle,\quad u\in\mathbb{R}^p.
\end{equation}
Here $\mathbb{D}: \mathbb{R}^p \rightarrow \mathbb{R}$ is a deterministic, non-negative, convex function with $\mathbb{D}(0) = 0$, $\mathbb{X}:\mathbb{R}^p\to\{0, \infty\}$ is the indicator function of some set $\Theta^\star\subseteq\mathbb{R}^p$, and $Y$ is a mean zero random variable in $\mathbb{R}^p$. Note that to understand the minimizers of $\mathbb{Z}(\cdot)$, we can without loss of generality assume convexity of $\mathbb{D}(\cdot)$. This is because the set of minimizers of any function $f:\mathbb{R}^p\to\overline{\mathbb{R}}$ is contained in the set of minimizers of the greatest convex minorant (GCM) of $f:\mathbb{R}^p\to\overline{\mathbb{R}}$. Formally, the greatest convex minorant of a function $f$ is defined as any function whose epigraph is the convex hull of the epigraph of the function $f$. Alternatively, for any $f:\mathbb{R}^p\to\overline{\mathbb{R}}$,
\begin{equation}\label{eq:def-GCM}
\mathrm{GCM}(f)(u) = \inf\left\{\sum_{j=1}^K \lambda_jf(u_j):\, K\ge 1,\, \sum_{j=1}^K \lambda_ju_j = u,\, \lambda_j \ge 0,\, \sum_{j=1}^K \lambda_j =1\right\}.
\end{equation}
See Proposition 2.31 of~\cite{rockafeller1998variational} for details. By definition, GCM of a convex function is itself and it is easy to verify (e.g., from~\eqref{eq:def-GCM}) that
\begin{equation}\label{eq:convexification-of-Z}
\mathrm{GCM}(\mathbb{Z})(u) = \widetilde{\mathbb{D}}(u) + \widetilde{\mathbb{X}}(u) + \langle u, Y\rangle,
\end{equation}
for some convex function $\widetilde{\mathbb{D}}(\cdot)$ and some convex indicator function $\widetilde{\mathbb{X}}(\cdot)$. In fact, $\tilde{\mathbb{X}}(\cdot)$ is the indicator function of the convex hull of $\Theta^\star$ and $\widetilde{\mathbb{D}}(u)$ is the greatest convex minorant of $\mathbb{D}(\cdot)$ as in~\eqref{eq:def-GCM}, except that $u_j$'s are restricted to $\Theta^\star$ (instead of $\mathbb{R}^p$).
This relation shows that one can always assume $\mathbb{D}(\cdot)$ and $\mathbb{X}(\cdot)$ are convex in~\eqref{eq:definition-Z}, to study the minimizers or the minimum value.

The following definitions are useful for our presentation. Let $P_Y$ denote probability measure induced by $Y$, i.e., $P_Y(A) = \mathbb{P}(Y \in A)$. For any function $f:\mathbb{R}^p\to\overline{\mathbb{R}}$, let $f^*(\cdot)$ denote the Legendre-Fenchel conjugate of $f$~\citep[Section 11.A]{rockafeller1998variational}, i.e.,
\[
f^*(y) = \sup_{x\in\mathbb{R}^p}\,\langle y, x\rangle - f(x).
\]
For any convex function $f:\mathbb{R}^p\to\overline{\mathbb{R}}$, let $\partial f(\cdot)$ denote the subdifferential mapping, i.e., 
\[
\partial f(u) := \{v\in\mathbb{R}^p:\, f(x) - f(u) \ge \langle v, x - u\rangle\,\mbox{for all}\, x\in\mathbb{R}^p\}.
\]
Any point in $\partial f(u)$ is called a subgradient of $f$ at $u$, and $v\in\partial f(u)$ if and only if $f(u) + f^*(v) = \langle u, v\rangle$. See, Section 2.4.3 of~\cite{bonnans2013perturbation} for details. 
Let $\mathcal{D}(f)$ denote the set of its differentiability points, i.e., $\mathcal{D}(f) = \{u\in\mathbb{R}^p:\, \partial f(u)\mbox{ is singleton}\}.$ A function $f:\mathbb{R}^p\to\overline{R}$ is called 
\begin{enumerate}
    \item {\em proper} if $f(x) < \infty$ for some $x\in\mathbb{R}^p$ and $f(x) > -\infty$ for all $x\in\mathbb{R}^p$;
    \item {\em lower semi-continuous} if for all $x\in\mathbb{R}^p$, $\liminf_{y\to x}f(y) = f(x)$;
    \item {\em weakly level-bounded} if for some $t \in \mathbb{R}$, $\{x:\, f(x) \le t\}$ is bounded and non-empty.
\end{enumerate}
See Chapter 1 of~\cite{rockafeller1998variational} for more details. Theorem 1.9 of~\cite{rockafeller1998variational} states that a function $f:\mathbb{R}^p\to\overline{\R}$ attains its infimum on $\mathbb{R}^p$ if $f(\cdot)$ is proper, lower semi-continuous, and weakly level bounded.

\subsection{Existence and Uniqueness of the Arg-Min}\label{sec:exist-arg-min}
Recall $\mathbb{Z}(u) = \mathbb{D}(u) + \mathbb{X}(u) + \langle u, Y\rangle$ for a non-negative, non-stochastic, convex function $\mathbb{D}:\mathbb{R}^p\to\mathbb{R}$ with $\mathbb{D}(0) = 0$, and a convex indicator function $\mathbb{X}(\cdot)$. With a slight abuse of notation, we assume $\Theta^\star$ is convex and that $\mathbb{X}(\cdot)$ is the indicator function of $\Theta^\star$. 

We will first discuss the existence of minimizes of $\mathbb{Z}(\cdot)$, i.e., the attainment of infimum of $\mathbb{Z}(\cdot)$ on $\mathbb{R}^p$. Note that the properties of $\mathbb{D}(\cdot)$ imply that it is a proper function. $\mathbb{X}(\cdot)$ is also a proper function if $\Theta^\star$ is non-empty, and $u\mapsto\langle u, Y\rangle$ is a proper function almost surely if $\mathbb{P}(\|Y\| = \infty) = 0$. Because $\mathbb{D}(\cdot)$ is a real-valued convex function, it is continuous on $\mathbb{R}^p$, and $\mathbb{X}(\cdot)$ is a lower-continuous function if and only if $\Theta^\star$ is closed; recall that all the tangent cones we discussed are closed. If $\Theta^\star$ is a bounded set, then $\mathbb{Z}(\cdot)$ is almost surely weakly level-bounded. If $\Theta^\star$ is unbounded, then level-boundedness requires some assumptions on the growth of $\mathbb{D}(\cdot)$ at infinity. In particular, we have the following result.
\begin{proposition}\label{prop:level-boundedness}
    If $\mathbb{D}(\cdot)$ is superlinear at $\infty$ in the sense that $\mathbb{D}(u)/\|u\| \to \infty$ as $\|u\| \to \infty$, then, with probability 1, $\mathbb{Z}(\cdot)$ is weakly level bounded. Moreover, if $\mathbb{D}(\cdot)$ is sublinear at $\infty$ in the sense that $\limsup_{\|u\|\to\infty}\mathbb{D}(u)/\|u\| \le \overline{c} < \infty$, and $\Theta^\star$ is unbounded in the sense that for some $v\in S^{p-1}$, $\Theta^\star \supset\{\lambda v:\, \lambda \ge 0\}$, then there exists a random vector $Y$ such that $\mathbb{Z}(\cdot)$ is not weakly level bounded with positive probability and consequently, the minimum is not attained.
\end{proposition}
See Section~\ref{appsec:proof-of-proposition-level-boundedness} for a proof. The proof shows that $Y$ can be taken as a mean zero Gaussian random vector for the second part of Proposition~\ref{prop:level-boundedness}. Here we mention that for a non-negative function, superlinearity at $\infty$ is same as coerciveness; see Defition 3.25 of~\cite{rockafeller1998variational}. Assuming level boundedness, Theorem 1.9 of~\cite{rockafeller1998variational} implies the following result. (The proof is omitted.)
\begin{theorem}\label{thm:existence}
    If $\mathbb{D}(\cdot)$ is non-negative convex function that is superlinear at $\infty$ and $\mathbb{X}(\cdot)$ is the indicator function of a closed convex set $\Theta^\star$, then for any random vector $Y$ with $\mathbb{P}(\|Y\| < \infty) = 1$, with probability 1, $\mathbb{Z}(\cdot)$ attains its infimum on $\mathbb{R}^p$.
\end{theorem}
It is important to note that Theorem~\ref{thm:existence} only provides sufficient conditions for existence of minimizers. From now on, we assume that $\mathbb{Z}(\cdot)$ attains its infimum.
The arg-min set of $u\mapsto\Z(u)$ over a set $\mathbb{R}^p$ is defined as
\begin{equation*}
    \mathcal{A}(\mathbb{Z}) := \{ s \in \mathbb{R}^p : \Z(s) = \inf_{u \in \mathbb{R}^p} \Z(u) \}.
\end{equation*}
The arg-min set is a singleton if the infimum is attained at a unique point. In such a case, we denote the argmin set by $W$. The following result provides a simple sufficient condition for the uniqueness of the minimizer. 
\begin{theorem} \label{thm:uniqueness}
    Assume that $\mathbb{Z}(\cdot)$ attains its infimum on $\mathbb{R}^p$. If $\mathbb{D}:\mathbb{R}^p\to\mathbb{R}$ is a strictly convex function on $\Theta^\star$ (i.e., $\mathbb{D}(\lambda u + (1-\lambda)v) < \lambda \mathbb{D}(u) + (1-\lambda)\mathbb{D}(v)$ for all $\lambda\in(0,1), u, v\in\Theta^\star$), then for any random vector $Y$, the argmin set $\mathcal{A}(\mathbb{Z})$ is a singleton. 
    
    Alternatively, if either $\Theta^\star$ is bounded or $u\mapsto \mathbb{D}(u)$ is superlinear at $\infty$, then for any absolutely continuous random vector $Y$, $\mathbb{Z}(\cdot)$ attains its infimum and the argmin set $\mathcal{A}(\mathbb{Z})$ is a singleton. 
\end{theorem}

See Section~\ref{appsec:proof-of-thm-uniqueness} for a proof. The second part of Theorem~\ref{thm:uniqueness} can also be derived from Theorem 11.8(d) of~\cite{rockafeller1998variational}. In the literature of convex analysis, there are results with weaker assumptions that imply uniqueness of the minimizer. We refer the interested reader to~\cite{kaufmann1988existence} and~\cite{planiden2016most,planiden2016strongly}. We believe for most statistical problems, Theorem~\ref{thm:uniqueness} suffices.
It is easy to see that in all class I examples we discussed, $\mathbb{D}(\cdot)$ is strictly convex and is superlinear: In Examples~\ref{ex:mean} and~\ref{ex:L2}, $\mathbb{D}(u) = \|u\|_2^2$ and $\mathbb{D}(u) = \langle u, Cu\rangle + \lambda_0\sum_{j=1}^p h_{\mu}(u^{(j)},\,\theta_0^{(j)})$. Because the quadratic functions are strictly convex and superlinear if the Hessian is positive definite, we get that the assumptions of Theorems~\ref{thm:existence} and~\ref{thm:uniqueness} are satisfied the minimum eigenvalue of $C$ is positive. In Example~\ref{ex:L1}, if, for example, $F_{i}(s) - F_i(0) = \lambda_i|s|^{\gamma}\mbox{sgn}(s)L(|s|)$ for a slowly varying function $L(\cdot)$ and $\gamma\in(0, \infty),$ then $\mathbb{D}(\cdot)$ is superlinear and strictly convex.

In the following subsections, we derive necessary and sufficient conditions for the symmetry of the distribution of the minimizer, assuming its existence and uniqueness. Recall that the argmin set $\mathcal{A}(\mathbb{Z})$ is denoted by $W$ under uniqueness. In Section~\ref{sec:1D}, we consider the one-dimensional case. In Section~\ref{sec:multi-dimensional-necessity}, we consider general multivariate case. This separation is useful for two reasons: (1) in the one-dimensional case, we can discuss median unbiasedness, in addition to symmetry (cf.~\eqref{eq:median-bias}); and (2) in the one-dimensional case, the derivation of necessary conditions is simplified because the minimizer can be written as a monotone function of $Y$, which is somewhat harder to use in the multivariate case. 

The following notation will be useful for our presentation. For any convex function $f:\mathbb{R}^p\to\mathbb{R}$, the directional derivative of $f$ at $x$ in the direction of $v\in S^{p-1}$ is defined as
\[
f'(x; v) := \lim_{t\downarrow 0}\, \frac{f(x + tv) - f(x)}{t}. 
\]
Convexity implies that such a limit always exists. It can be verified that $f'(x; v) = \sup_{g\in\partial f(x)}\,g^{\top}v$, where $\partial f(x)$ represents the subdifferential of $f$ at $x$. If $p = 1$, there are only two directions $v = -1$ or $v = +1$ and we write the left and right derivatives of $f$ at $x$ as
\[
f_+'(x) = \lim_{t\downarrow 0}\, \frac{f(x + t) - f(x)}{t},\quad\mbox{and}\quad f_-'(x) = \lim_{t\uparrow 0}\,\frac{f(x + t) - f(x)}{t}.
\]
Note that $f_+'(x) = f'(x; +1)$ and $f_-'(x) = -f'(x; -1)$.

\subsection{One-dimensional Case} \label{sec:1D}
The following are the main results for the one-dimensional case. Recall that we are assuming $\mathbb{D}(\cdot)$ is non-negative convex function with $\mathbb{D}(0) = 0$, and $\mathbb{X}(\cdot)$ is the indicator function for a convex set $\Theta^\star$ (which in the one-dimensional case is an interval).
\begin{theorem}\label{thm:1d-sym-A1-necessity}
    Assume that $W$ is well-defined. If $Y$ is a symmetric random variable whose measure dominates the Lebesgue measure on $\mathbb{R}$, then $W$ is symmetrically distributed around zero if and only if~\ref{assump:even-ness-sufficiency} holds true. 
\end{theorem}

\begin{theorem} \label{thm:1d-med} 
    If $W$ is well-defined, $0\in\Theta^\star$, and $Y$ is median unbiased, then $W$ is median unbiased.

    Moreover, if $W$ is well-defined, $0$ is in the interior of $\Theta^\star$, and $W$ is median unbiased, then $Y$ is median unbiased.
\end{theorem}
The proofs of Theorems~\ref{thm:1d-sym-A1-necessity} and~\ref{thm:1d-med} can be found in Sections~\ref{appsec:proof-of-1d-sym-A1-necessity} and~\ref{appsec:proof-of-1d-med}, respectively. 
Recall that for the general stochastic process $\mathbb{Z}(\cdot)$ (defined in~\eqref{eq:definition-Z}) with a potentially non-convex $\mathbb{D}(\cdot)$ and $\mathbb{X}(\cdot)$, the results above continue to hold true if $\mathbb{D}(\cdot)$ and $\mathbb{X}(\cdot)$ are replaced with their convexifications $\widetilde{\mathbb{D}}(\cdot)$ and $\widetilde{\mathbb{X}}(\cdot)$; see the discussion following~\eqref{eq:convexification-of-Z}. The assumption that $Y$ is supported on all of $\mathbb{R}$ can be relaxed. From the proof, it follows that if $W$ is symmetrically distributed, then the conjugates of $u\mapsto \D(u) + \X(u)$ and $u\mapsto \D(-u) + \X(-u)$ much match on the support of $Y$. What this implies about assumption~\ref{assump:even-ness-sufficiency} is unclear, at present. 

Theorem~\ref{thm:1d-med} proves that median unbiased holds under much weaker conditions than~\ref{assump:even-ness-sufficiency}. Although we haven't explored median unbiasedness in the multidimensional case, Theorem~\ref{thm:1d-med} can still be useful in the multivariate case, if $\mathbb{D}(\cdot)$ and $\mathbb{X}(\cdot)$ are separable. For instance, in Example~\ref{ex:mean}, $\mathbb{D}(u) = \sum_{j=1}^p \{(u^{(j)})^2 - 2u^{(j)}Y^{(j)}\}$ and if $\Theta$ is a separable constraint set, i.e., $\Theta = \{u\in\mathbb{R}^p:\, u^{(j)}\in \Theta_j \subseteq \mathbb{R}\}$, then symmetry and median unbiasedness of $W$ can be understood by applying Theorem~\ref{thm:1d-sym-A1-necessity},~\ref{thm:1d-med}, because the minimizer can be obtained by considering $p$ univariate minimization problems. The same strategy works in Example~\ref{ex:L2} if $C$ is a diagonal matrix.

Regarding the assumptions on $Y$ in Theorem~\ref{thm:1d-sym-A1-necessity}, we note that in class I limiting distributions $Y$ arises as the limiting distribution of properly scaled $\nabla \mathbb{M}_n(\theta_0) - \mathbb{E}[\nabla\mathbb{M}_n(\theta_0)]$. Hence, $Y$ will be an infinitely divisible distributions. It is well-known that infinitely divisible distributions cannot be supported on a bounded set and any absolutely continuous infinitely divisible distribution with support of $\mathbb{R}$ dominates the Lebesgue measure on $\mathbb{R}$.
\subsection{Multidimensional Case}\label{sec:multi-dimensional-necessity}

We now begin our analysis of the general set up. Firstly, a lemma concerning the extraction of a measurable selection of the subdifferential mapping is presented.

\begin{lemma} \label{lem:selection}
    If $f : \R^p \rightarrow \R \cup \{+ \infty \}$ is a convex function such that $\partial f$ is a surjective map, then there exists a measurable selection of the subdifferential mapping of its conjugate, i.e., a borel function $\widetilde{\nabla}f^* : \R^p \rightarrow \R^p$ such that $\widetilde{\nabla} f^* (u) \in \partial f^* (u)$ for all $u \in \R^p$.
\end{lemma} 

Next, we exhibit a result which, under mild assumptions on the random shift, cements the necessity of $d$ being an even function for the symmetry of the arg-min distribution. Note that these assumptions hold true when $Y$ is a mean zero Gaussian.

\begin{theorem} \label{thm:gen-sym}
    Assume that $\D$ is superlinear, $\Theta^\star = \mathbb{R}^p$ and $Y$ is an absolutely continuous, centrally (resp. spherically, sign) symmetric random vector that dominates the Lebesgue measure on $\mathbb{R}^p$. Then $W$ is centrally (resp. spherically, sign) symmetric if and only if~\ref{assump:even-ness-sufficiency} holds true.
\end{theorem}

The proofs of Theorem~\ref{thm:1d-sym-A1-necessity} and~\ref{thm:gen-sym} make use of the fact that there is a unique optimal transport map that transports an absolutely continuous random variable/vector to a fixed distribution. This is an unexpected application of optimal transport theory in the study of symmetry of limiting distributions.

\section{Conclusions}\label{sec:conclusions}

In this paper, we have described three commonly occurring classes of limiting distributions of $M$-estimators and provided conditions under which these limiting distributions are symmetric. We further examined the necessity of these conditions for symmetry in one simple class of limiting distributions. The classical results of optimal transport theory makes appearance in our study, which is unexpected. 

More research in this direction would significantly aid statistical inference for some of the irregular or non-standard statistical problems. For example, it would be interesting to find necessary and sufficient conditions for the minimizer of a drifted Brownian motion to be symmetric. In addition to symmetry, other properties of limiting distributions such as median unbiasedness and/or unimodal can also yield simple statistical inference procedures~\citep{kuchibhotla2024hulc}. It is not clear how to transfer some of the techniques in this paper to study conditions for other properties of minimizers of stochastic processes. 

\bibliographystyle{apalike}
\bibliography{references}
\newpage

\setcounter{section}{0}
\setcounter{equation}{0}
\setcounter{figure}{0}
\renewcommand{\thesection}{S.\arabic{section}}
\renewcommand{\theequation}{E.\arabic{equation}}
\renewcommand{\thefigure}{A.\arabic{figure}}
\renewcommand{\theHsection}{S.\arabic{section}}
\renewcommand{\theHequation}{E.\arabic{equation}}
\renewcommand{\theHfigure}{A.\arabic{figure}}
  \begin{center}
  \Large {\bf Supplement to ``On the Symmetry of Limiting Distribution of M-estimators''}
  \end{center}
       
\begin{abstract}
This supplement contains the proofs of all the main results in the paper, and some supporting lemmas. 
\end{abstract}

\section{Auxiliary Results}\label{appsec:proof-of-median-unbiasedness}

\begin{claim}
    If $r_n(\widehat{\theta}_n - \theta(P)) \overset{d}{\to} W_P$ for $P\in\mathcal{P}$, then 
    \begin{equation*} 
        \liminf_{n \rightarrow \infty} \, \medbias_P(\widehat{\theta}_n; \theta(P)) \geq
        \left(\frac{1}{2} - \min_{s\in\{-1, 1\}} \mathbb{P}_P (sW_P \ge 0) \right)_+ . 
    \end{equation*}
    Further, if $\mathbb{P}_P(W_P = 0) = 0$, then limit of $\medbias_P(\widehat{\theta}_n; \theta(P))$ exists and 
    \begin{equation*}
         \lim_{n \rightarrow \infty} \, \medbias_P(\widehat{\theta}_n; \theta(P)) =
        \left(\frac{1}{2} - \min_{s\in\{-1, 1\}} \mathbb{P}_P (sW_P \ge 0) \right)_+ .
    \end{equation*}
\end{claim}

\begin{proof}
    Since $r_n(\widehat{\theta}_n - \theta(P)) \overset{d}{\to} W_P$ for $P\in\mathcal{P}$, Portmanteau theorem implies that
    \begin{align*}
        &\limsup_{n \rightarrow \infty}\, \P_P(r_n(\widehat{\theta}_n - \theta(P)) \geq 0) \leq \P_P (W_P \geq 0) \\
        &\limsup_{n \rightarrow \infty}\, \P_P(r_n(\widehat{\theta}_n - \theta(P)) \leq 0) \leq \P_P (W_P \leq 0).
    \end{align*}
    Consequently, we have 
    \begin{align*}
        \min_{s\in\{-1, 1\}} \mathbb{P}_P (sW_P \ge 0) &\geq 
        \min_{s\in\{-1, 1\}} \limsup_{n \rightarrow \infty}\, \P_P(s(\widehat{\theta}_n - \theta(P)) \geq 0) \\
        &\geq  \limsup_{n \rightarrow \infty}\, \min_{s\in\{-1, 1\}}
        \P_P(s(\widehat{\theta}_n - \theta(P)) \geq 0),
    \end{align*}
    which implies that
    \begin{align*}
        \left(\frac{1}{2} - \min_{s\in\{-1, 1\}} \mathbb{P}_P (sW_P \ge 0) \right)_+ &\leq 
        \left(\liminf_{n \rightarrow \infty} \left( \frac{1}{2} - \min_{s\in\{-1, 1\}} \P_P(s(\widehat{\theta}_n - \theta(P)) \geq 0) \right) \right)_+ \\
        &\leq \: \liminf_{n \rightarrow \infty} \left( \frac{1}{2} - \min_{s\in\{-1, 1\}} \P_P(s(\widehat{\theta}_n - \theta(P)) \geq 0) \right)_+,
    \end{align*}
    completing the first part of the proof. Next, if we know that $\mathbb{P}_P(W_P = 0) = 0$, then 0 is a continuity point of the distribution of $W_P$ under $P$ and hence
    \begin{align*}
        &\P_P(r_n(\widehat{\theta}_n - \theta(P)) \geq 0) \rightarrow \P_P (W_P \geq 0) \\
        &\P_P(r_n(\widehat{\theta}_n - \theta(P)) \leq 0) \rightarrow \P_P (W_P \leq 0),
    \end{align*}
    which immediately gives
    \begin{equation*}
        \left( \frac{1}{2} - \min_{s\in\{-1, 1\}} \P_P(s(\widehat{\theta}_n - \theta(P)) \geq 0) \right)_+ \! \longrightarrow 
        \left(\frac{1}{2} - \min_{s\in\{-1, 1\}} \mathbb{P}_P (sW_P \ge 0) \right)_+ ,
    \end{equation*}
    as desired.
\end{proof}

\section{Proofs of results in \Cref{sec:sufficiency}}

\subsection{Proof of \Cref{suff:NS}}

If $S(u) = \langle u, Y \rangle$ is even, then $u \mapsto \langle u, Y \rangle$ and $u \mapsto \langle -u, Y \rangle$ have the same law. In particular, we must have $\langle u, Y \rangle \ed \: \langle -u, Y  \rangle = \langle u, -Y  \rangle$ for all $u \in \R^p$. Consequently, Cramer-Wold theorem implies that $Y \ed -Y$. To prove the converse, assume that $Y$ is symmetric and consider the random vector $(\mathbb{S}(u_1), \hdots, \mathbb{S}(u_k))$ for any $k \geq 1$ and $u_1, \hdots, u_k \in \R^p$. It suffices to show that this has the same distribution as $(\mathbb{S}(-u_1), \hdots, \mathbb{S}(-u_k))$. So we write for $t_1, \hdots, t_k \in \R$,
\begin{align*}
    \E \left[\exp \left(i \sum_{j=1}^{k} t_j \mathbb{S}(u_j) \right) \right] &=  \E \left[\exp \left(i \sum_{j=1}^{k} t_j\langle u_j, Y \rangle \right) \right]   \\
    &=  \E \left[\exp \left(i \sum_{j=1}^{k} t_j\langle u_j, -Y \rangle \right) \right]   \tag{from symmetry of Y}\\
    &=  \E \left[\exp \left(i \sum_{j=1}^{k} t_j\langle -u_j, Y \rangle \right) \right]   \\
    &= \E \left[\exp \left(i \sum_{j=1}^{k} t_j \mathbb{S}(-u_j) \right) \right] .
\end{align*}
This shows that the characteristic functions of the two random vectors coincide everywhere, which implies that they are identically distributed, as required. 
    
\subsection{Proof of \Cref{suff:GP}}

By definition, the Gaussian process $\mathbb{S}(u) = \mathrm{GP}_{\Sigma}(u)$ is even if $u \mapsto \mathrm{GP}_{\Sigma}(u)$ and $u \mapsto \mathrm{GP}_{\Sigma}(-u)$ share the same law. Since both of these are mean zero Gaussian processes, this is equivalent to them having identical covariance operators. The map $u \mapsto \mathrm{GP}_{\Sigma}(-u)$ has covariance operator $\Sigma'(u, v) = \Sigma(-u, -v)$, and hence $\mathbb{S}(u)$ is even if and only if $\Sigma(u, v) = \Sigma(-u, -v)$ for all $u, v \in \R^p$. 

\subsection{Proof of \Cref{suff:GPHP}}

Suppose that $\nu(A, B) = \nu (A, -B)$ for all Borel measurable sets $A\subseteq \R, B \subseteq \R^p$. By definition of the Poisson hyperplane process, we have
\begin{equation*}
    \mbox{Poisson}(\gamma, \nu)(u) = \sum_{i=1}^{\infty} V_i \mathbf{1}\{ \tau_i \leq \langle u, U_i\rangle\}, \qquad u \in \R^p,
\end{equation*}
where $(V_i, U_i)_{i \geq 1}$ is an independent sequence of random vectors distributed according to $\nu$. For $k \geq 1$ and $u_1, \hdots, u_k \in \R^p$, we can view $(\mbox{Poisson}(\gamma, \nu)(u_1), \hdots, \mbox{Poisson}(\gamma, \nu)(u_k))$ as the output of the measurable map
\begin{equation*}
    (V_i, U_i)_{i \geq 1} \mapsto \left(\sum_{i=1}^{\infty} V_i \mathbf{1}\{ \tau_i \leq \langle u_1, U_i\rangle\} , \hdots, \sum_{i=1}^{\infty} V_i \mathbf{1}\{ \tau_i \leq \langle u_k, U_i\rangle\}\right).
\end{equation*}
The assumption on $\nu$ implies that 
\begin{equation*}
    (V_i, U_i)_{i \geq 1} \ed \; (V_i, -U_i)_{i \geq 1}.
\end{equation*}
Since $(V_i, U_i)_{i \geq 1}$ is independent of $(\tau_i)_{i \geq 1}$, we can then condition on the $\tau_i$'s to obtain
\begin{equation*}
    (\mbox{Poisson}(\gamma, \nu)(u_1), \hdots, \mbox{Poisson}(\gamma, \nu)(u_k)) \ed \; (\mbox{Poisson}(\gamma, \nu)(-u_1), \hdots, \mbox{Poisson}(\gamma, \nu)(-u_k)),
\end{equation*}
which proves that the stochastic process $\mbox{Poisson}(\gamma, \nu)(u)$ is even. This trivially implies that 
\begin{equation*}
    \E[\mbox{Poisson}(\gamma, \nu)(u)] = \E[\mbox{Poisson}(\gamma, \nu)(-u)],
\end{equation*}
from which we can conclude that the centered Poisson hyperplane process $\mathbb{S}(u) = \mbox{cPoisson}(\gamma, \nu)(u)$ is even under the given assumption.

\section{Proofs of results in \Cref{sec:necessity}}
\subsection{Proof of Proposition~\ref{prop:level-boundedness}}\label{appsec:proof-of-proposition-level-boundedness}
Observe that for any $t\in\mathbb{R}$,
\begin{align*}
\{u\in\mathbb{R}^p:\, \mathbb{Z}(u) \le t\} ~\subseteq~ \{u\in\Theta^\star:\,\|u\|(\mathbb{D}(u)/\|u\| - \|Y\|) \le t\}.
\end{align*}
If $\mathbb{D}(u)$ is superlinear, then with probability 1,
\[
\lim_{\|u\|\to\infty}\, \frac{\mathbb{D}(u)}{\|u\|} - \|Y\| = \infty,
\]
which implies that $\|u\|(\mathbb{D}(u)/\|u\| - \|Y\|)$ cannot stay bounded by $t$ unless $\|u\|$ is bounded. Hence, $\mathbb{Z}(\cdot)$ is weakly level-bounded.

Note that for any $t\in\mathbb{R}$,
\[
\{u\in\mathbb{R}^p:\,\mathbb{Z}(u) \le t\} \supseteq \{\lambda v:\, \lambda(\mathbb{D}(\lambda v)/\lambda + \langle v, Y\rangle) \le t\}.
\]
Take any random vector $Y$ such that $\mathbb{P}(\langle v, Y\rangle \le -\overline{c} - 1) \ge \delta$ for some $\delta > 0$; such random vectors exist as one can take $Y$ to be a mean zero Gaussian with $(\overline{c}+1)^2I_p$ as the covariance matrix. Also, because $\mathbb{D}(\cdot)$ is sublinear, there exists $\lambda^*>0$ such that for all $\lambda \ge \lambda^*$, $\mathbb{D}(\lambda v) \le \lambda\overline{c}$. Hence, with probability at least $\delta$, 
\[
\lambda(\mathbb{D}(\lambda v)/\lambda + \langle v, Y\rangle) \le -\lambda\quad\mbox{for all}\quad \lambda \ge \lambda^*,
\]
which implies that with probability at least $\delta$, $\{\lambda v:\, \lambda \ge \max\{\lambda^*, -t\}\}$ (an unbounded set) belongs to $\{u:\mathbb{Z}(u) \le t\}$. Therefore, $\mathbb{Z}(\cdot)$ is not weakly level bounded. From the analysis above, it is also clear that
\[
\min_{u\in\mathbb{R}^p} \mathbb{Z}(u) \le \min_{\lambda \ge \lambda^*}\, \mathbb{D}(\lambda v) + \lambda \langle v, Y\rangle \overset{(a)}{\le} \min_{\lambda \ge \lambda^*} -\lambda = -\infty, 
\]
where inequality (a) above holds with probability at least $\delta$. Given that $|\mathbb{Z}(u)| < \infty$ for any $u\in\mathbb{R}^p$, the minimum is not attained with probability at least $\delta$.

\subsection{Proof of \Cref{thm:uniqueness}}\label{appsec:proof-of-thm-uniqueness}

The proof follows by noting that strict convexity of $\mathbb{D}(\cdot)$ implies strict convexity of $u\mapsto\mathbb{D}(u) + \langle u, Y\rangle$. Minimizer of a strict convex function on a convex set is unique, if it exists. This is a standard result in convex analysis. We provide this short proof for completeness. 

Suppose that $\Z$ attains its infimum on $\Theta^\star$, and $u, v \in \Theta^\star$ be two distinct minimizers, i.e.,
\begin{equation*}
    \Z(u) = \Z(v) = \inf_{s \in \Theta^\star} \Z(s).
\end{equation*}
Since $\Theta^\star$ is a convex set, $\lambda u + (1-\lambda)v \in \Theta^\star$ for any $\lambda \in (0,1)$. By strict convexity of $\Z$ on $\Theta^*$, we have 
\begin{equation*}
    \Z(\lambda u + (1-\lambda)v) < \lambda \Z(u) + (1-\lambda) \Z(v) = \inf_{s \in \Theta^\star} \Z(s),
\end{equation*}
which is clearly impossible. Thus the minimizer, if it exists, is unique.

    Since $\Z$ takes value $+\infty$ outside $\Theta^\star$, minimizing $\Z$ over $\Theta^\star$ is equivalent to minimizing it over $\R^p$. Note that $\Z$ is minimized at a unique point in $\R^p$ if and only if $\Z^\star$ is differentiable at 0 \citep[Theorem 27.1]{Rockafellar+1970}. Also,
\begin{align*}
    \Z^\star(u) &= \sup_{x \in \R^p} \; \langle u, x \rangle - \Z(x) \\
    &= \sup_{x \in \R^p} \; \langle u, x \rangle - \D(x) - \X(x) - \langle x, Y \rangle \\
    &= \sup_{x \in \R^p} \; \langle u-Y, x \rangle - \D(x) - \X(x) \\
    &= (\D + \X)^\star (u-Y).
\end{align*}
Hence, differentiability of $\Z^*$ at $0$ is same as differentiability of the convex conjugate of $(\mathbb{D} + \mathbb{X})$ at $Y$. From Theorem 11.8 of~\cite{rockafeller1998variational}, it follows that the domain of $(\D + \X)^*(\cdot)$ is $\mathbb{R}^p$ if $u\mapsto \mathbb{D}(u) + \mathbb{X}(u)$ is coercive or superlinear, which is our assumption. Because convex functions are almost everywhere differentiable on their doamin and $v\mapsto (\mathbb{D} + \X)^*(u)$ is a convex function, the assumption of absolute continuity of $Y$ implies that with probability 1, $\mathbb{Z}^*(\cdot)$ is differentiable at $0$.

\subsection{Proof of \Cref{thm:1d-sym-A1-necessity}}\label{appsec:proof-of-1d-sym-A1-necessity}
    By Theorem~11.8(b) of~\cite{rockafeller1998variational}, we get $W$ (the unique minimizer of $\mathbb{Z}(\cdot)$) is the derivative of the conjugate of $u\mapsto \mathbb{D}(u) + \mathbb{X}(u) + \langle u, Y\rangle$ at zero. Similarly, $-W$, the unique minimizer of $u\mapsto \mathbb{Z}(-u)$, is the derivative of the conjugate of $u\mapsto \mathbb{D}(-u) + \mathbb{X}(-u) - \langle u, Y\rangle$. For notational convenience, set
    \[
    g(v) = \sup_{u\in\mathbb{R}^p} \langle u, v\rangle - \mathbb{D}(u) - \mathbb{X}(u),\quad\mbox{and}\quad h(v) = \sup_{u\in\mathbb{R}^p}\, \langle u, v\rangle - \D(-u) - \X(-u),
    \]
    as the conjugates of $u\mapsto \mathbb{D}(u) + \mathbb{X}(u)$ and $u\mapsto\mathbb{D}(-u) + \mathbb{X}(-u)$.
    The uniqueness of the minimizers imply that
    Then $W = g'(-Y)$ and $-W = h'(Y)$; the derivatives here exist almost surely from almost sure uniqueness of minimizers. This implies that $g'(-Y) \overset{d}{=} h'(Y)$. Because of the symmetry of $Y$, this is equivalent to $g'(-Y) \overset{d}{=} h'(-Y)$. Because $u\mapsto g'(-u)$ and $u\mapsto h'(-u)$ are non-increasing functions and they both yield the same distribution when evaluated at $Y$, they must be equal almost surely on the support of $Y$, if $Y$ is absolutely continuous with respect to the Lebesgue measure on its support. This, for example, follows from the uniqueness of the optimal transport map in the one-dimensional case; see Theorem 3.1 of~\cite{Ambrosio2003}. Because the support of $Y$ is $\mathbb{R}^p$, and equality of convex conjugates implies equality of functions~\citep[Theorem 11.1]{rockafeller1998variational}, the result follows. Note that, in general, $W$ is symmetrically distributed if and only if the conjugates $g$ and $h$ are equal on the support of $Y$.

\subsection{Proof of \Cref{thm:1d-med}}\label{appsec:proof-of-1d-med}
By the first order necessary condition for minimizers, 
\[
(\D + \X)_-'(W) + Y \le 0 \le (\mathbb{D} + \mathbb{X})_+'(W) + Y.
\]
Because the left and right derivatives are monotone non-decreasing, 
\[
\{W > 0\} \subseteq \{Y < -(\D + \X)_+'(0)\}\quad\mbox{and}\quad \{W < 0\} \subseteq \{Y > -(\D + \X)_-'(0)\}.
\]
Because $\Theta^\star$ is a closed convex set, let $\Theta^\star$ be $[a, b]$. It is easy to verify that
\[
(\D + \X)_+'(0) = \begin{cases}\D_+'(0), &\mbox{if }0 \le a < b,\\
\infty, &\mbox{if }a < 0 = b,\end{cases}\quad\mbox{and}\quad (\D + \X)_-'(0) = \begin{cases}
-\infty, &\mbox{if }0 = a < b,\\
\D_-'(0), &\mbox{if }a < 0 \le b,\\
\end{cases}
\]

Now, recall that $\D(\cdot)$ attains its minimum at zero and hence $\D_+'(0) \ge 0\ge \D_-'(0)$. Therefore,
\[
\mathbb{P}(W > 0) \le \mathbb{P}(Y < -(\D + \X)_+'(0)) \le \begin{cases}
\mathbb{P}(Y < -\D_+'(0)) \le \mathbb{P}(Y < 0), &\mbox{if }0\le a < b,\\
\mathbb{P}(Y < -\infty) = 0, &\mbox{if } 0 = b.
\end{cases}
\]
Similarly, 
\[
\mathbb{P}(W < 0) \le \mathbb{P}(Y > -(\D + \X)_-'(0)) \le \begin{cases}
\mathbb{P}(Y > -\D_-'(0)) \le \mathbb{P}(Y > 0), &\mbox{if }a < 0 \le b,\\
\mathbb{P}(Y > \infty) = 0, &\mbox{if }0 = a.
\end{cases}
\]
Therefore, 
\[
\max\{\mathbb{P}(W < 0),\, \mathbb{P}(W > 0)\} \le \max\{\mathbb{P}(Y < 0),\,\mathbb{P}(Y > 0)\}.
\]
Median unbiasedness of $Y$ implies that the right hand side is at most $1/2$ and hence, $W$ is also median unbiased.

On the other hand, observe that if $0 \in (a, b)$, then
\begin{align*}
\{W \ge 0\} &\subseteq \{Y \le -(\D + \X)_+'(0)\} = \{Y \le -\D_+'(0)\} \subseteq \{Y \le 0\},\\
\{W \le 0\} &\subseteq \{Y \ge -(\D + \X)_-'(0)\} = \{Y \ge -\D_-'(0)\} \subseteq \{Y \ge 0\}.
\end{align*}
Here the final inclusions follow from $\D_-'(0) \le 0 \le \D_+'(0)$. Therefore, 
\[
\min\{\mathbb{P}(W \le 0),\,\mathbb{P}(W \ge 0)\} \le \min\{\mathbb{P}(Y \le 0),\,\mathbb{P}(Y \ge 0)\}.
\]
This shows that median unbiasedness of $W$ implies median unbiasedness of $Y$.



\subsection{Proof of \Cref{lem:selection}}

Since $\partial f$ is surjective, $\partial f^* (u)$ is non-empty for all $u \in \R^p$. This implies that the effective domain of $f^*$ is whole of $\R^p$, that is, $f^*$ is real valued. As a consequence, $f^*$ is continuous everywhere and hence, $\partial f^*$ is upper-semicontinuous as a correspondence. This implies that it is also weakly measurable, and by the selection theorem of Kuratowski and Ryll-Nardzewski \cite[Theorem~5.1]{Parthasarathy1972}, there exists a borel measurable selection $\widetilde{\nabla}f^*$. \\

\subsection{Proof of \Cref{thm:gen-sym}}

While we will prove the theorem assuming $Y$ is centrally symmetric, the other scenarios will follow the same idea. Firstly, note that $W$ is well defined by~\Cref{thm:uniqueness}. The if part follows from~\ref{thm:suff}.

The converse will be proved following the technique of Theorem~\ref{thm:1d-sym-A1-necessity}, by utilizing results about optimal transport. To begin, note that $\d$ is superlinear and hence $\partial \d$ is surjective. By \Cref{lem:selection}, we can extract a measurable selection $\widetilde{\nabla}\d^*$ of $\partial \d^*$. Observe that when $u \in \mathcal{D}(\d^*)$, $\widetilde{\nabla}\d^*(u) = \nabla \d^* (u)$. For any Borel set $A$, 
    \begin{align}
        \P \big(\widetilde{\nabla}\d^* (Y) \in A \big) &= \P \big(\widetilde{\nabla}\d^* (-Y) \in A \big) \label{step1}\\
        &= \P \big(\widetilde{\nabla}\d^* (-Y) \in A, -Y \in D(\d^*) \big) + \P \big(\widetilde{\nabla}\d^* (-Y) \in A, -Y \notin D(\d^*) \big) \notag \\
        &= \P \big(\nabla \d^* (-Y) \in A, -Y \in D(\d^*) \big) + \P \big(\widetilde{\nabla}\d^* (-Y) \in A, -Y \notin D(\d^*) \big) \notag \\
        &= \P \big(W \in A, -Y \in D(\d^*) \big) + \P \big(\widetilde{\nabla}\d^* (-Y) \in A, -Y \notin D(\d^*) \big) \label{step2}\\
        &= \P(W \in A). \label{step3}
    \end{align}
Here, $\eqref{step1}$ is due to central symmetry of $Y$, $\eqref{step2}$ is due to the fact that $W = \nabla \d^* (-Y)$ w.p. 1 and \eqref{step3} follows from the fact that $Y$ is absolutely continuous and the set of non-differentiability points of $\d^*$ has Lebesgue measure zero. By repeating this calculation again, one can show 
    \begin{align*}
        \P \big(-\widetilde{\nabla}\d^* (-Y) \in A \big) &= \P(-W \in A) \\
        &= \P(W \in A)     \\
        &= \P \big(\widetilde{\nabla}\d^* (Y) \in A \big),
    \end{align*}
where the second last step is due to central symmetry of $W$. Thus, we obtained
    \begin{equation*}
        W \ed \: \widetilde{\nabla}\d^* (Y) \ed \: -\widetilde{\nabla}\d^* (-Y).
    \end{equation*}

If we define $\g(u) := \d(-u)$, then we have $\g^*(u) = \d^*(-u)$ and consequently, $\partial \g^* (-u) = -\partial \d^*(-u)$ as sets. Thus, we can obtain a measurable selection $\widetilde{\nabla}\g^*$ of $\partial \g^*$ such that $\widetilde{\nabla}\g^*(u) = -\widetilde{\nabla}\d^* (-u)$ for all $u$. So, we get
    \begin{equation*}
        W \ed \: \widetilde{\nabla}\d^* (Y) \ed \: \widetilde{\nabla}\g^* (Y).
    \end{equation*}
Now, $\widetilde{\nabla}\d^*$ and $\widetilde{\nabla}\g^*$ are selections of the subdifferentials of convex functions, which transport an absolutely continuous measure $P_Y$ onto $P_W$. From Briener-McCann theorem \cite[Theroem~2.32]{villani2021topics} in optimal transport, we conclude that $\widetilde{\nabla}\d^*$ and $\widetilde{\nabla}\g^*$ coincide $P_Y$ almost everywhere. $Y$ has support $\R^p$ implies that the set of points where they differ is at most a set of Lebesgue measure zero.
    
From the mean value theorem for convex functions \cite[Theorem~2.3.4]{hiriart2004fundamentals}, we have
    \begin{equation*}
        \d^*(u) = \int_0^1 \langle \, \widetilde{\nabla}\d^* (yu), u \, \rangle \: dy \quad \text{ and } \quad \g^*(u) = \int_0^1 \langle \, \widetilde{\nabla}\g^* (yu), u \, \rangle \: dy.
    \end{equation*}
The argument in the last paragraph shows that $\d^*(u) = \g^*(u)$ everywhere. Consequently, $\d(u) = \g(u)$ everywhere, thus proving that $\d$ is even.   

\end{document}